\documentclass[aap,MSNbibl,citesort,dvips]{arximspdf}
\usepackage{mathbh}
\usepackage{graphicx}

\doi{10.1214/11-AAP796}
\volume{22}
\issue{4}
\pubyear{2012}
\firstpage{1362}
\lastpage{1410}

\makeatletter
\newcommand{\rref}[1]{\mbox{{\fontsize{8.36}{10}\selectfont{\ref{#1}}}}}

\newcommand{\N}{\mathbb N}
\newcommand{\Z}{\mathbb{Z}}
\newcommand{\Zd}{\mathbb{Z}^d}
\newcommand{\Q}{\mathbb{Q}}
\newcommand{\R}{\mathbb{R}}
\newcommand{\Rd}{\mathbb{R}^d}
\renewcommand{\P}{\mathbb{P}}
\newcommand{\E}{\mathbb{E}}
\newcommand{\Ed}{\mathbb{E}^d}
\newcommand{\Card}[1]{\vert #1 \vert}
\newcommand{\one}{\mathbh{1}}

\newcommand{\Ebarre}{\overline{\mathbb{E}}}
\newcommand{\Pbarre}{\overline{\mathbb{P}}}
\newcommand{\Qbarre}{\overline{\mathbb{Q}}}
\newcommand{\D}{\mathcal{D}}
\newcommand{\T}[2]{{#1}.{#2}}

\renewcommand{\epsilon}{\varepsilon}
\renewcommand{\limsup}{\mathop{\overline{\lim}}}
\renewcommand{\liminf}{\mathop{\underline{\lim}}}

\newcommand{\Erg}{\operatorname{Erg}}

\newtheorem{theorem}{Theorem}
\newtheorem{lemme}[theorem]{Lemma}
\newtheorem{coro}[theorem]{Corollary}
\newtheorem{prop}[theorem]{Proposition}

\newproclaim{Remark}{Remark}

\makeatother

\begin{document}
\begin{frontmatter}

\title{Asymptotic shape for the contact process in random environment}
\runtitle{Contact process in random environment}

\begin{aug}
\author[A]{\fnms{Olivier} \snm{Garet}\corref{}\ead[label=e1]{Olivier.Garet@iecn.u-nancy.fr}\ead[label=u1,url]{http://www.iecn.u-nancy.fr/\textasciitilde garet/}}
\and
\author[A]{\fnms{R\'{e}gine} \snm{Marchand}\ead[label=e2]{Regine.Marchand@iecn.u-nancy.fr}\ead[label=u2,url]{http://www.iecn.u-nancy.fr/\textasciitilde rmarchan/}}
\runauthor{O. Garet and R. Marchand}
\affiliation{University of Lorraine}
\address[A]{Institut \'Elie Cartan (Math{\'e}matiques)\\
Universit{\'e} de Lorraine\\
Campus Scientifique, BP 239 \\
54506 Vandoeuvre-l{\`e}s-Nancy Cedex\\
France\\
\printead{e1}\\
\hphantom{E-mail: }\printead*{e2}\\
\printead{u1}\\
\hphantom{URL: }\printead*{u2}} 
\end{aug}

\received{\smonth{5} \syear{2010}}
\revised{\smonth{3} \syear{2011}}

%
\begin{abstract}
The aim of this article is to prove asymptotic shape theorems for the
contact process in stationary random environment. These theorems
generalize known results for the classical contact process. In
particular, if $H_t$ denotes the set of already occupied sites at time
$t$, we show that for almost every environment, when the contact
process survives, the set $H_t/t$ almost surely converges to a compact
set that only depends on the law of the environment. To this aim, we
prove a~new almost subadditive ergodic theorem.
\end{abstract}

%
\begin{keyword}[class=AMS]
\kwd[Primary ]{60K35}
\kwd[; secondary ]{82B43}.
\end{keyword}
\begin{keyword}
\kwd{Random growth}
\kwd{contact process}
\kwd{random environment}
\kwd{almost subadditive ergodic theorem}
\kwd{asymptotic shape theorem}.
\end{keyword}

\end{frontmatter}

\section{Introduction}\label{intro}
The aim of this paper is to obtain an asymptotic shape theorem for the
contact process in random environment on $\Zd$. The ordinary contact
process is a famous interacting particle system modeling the spread of
an infection on the sites of $\Zd$. In the classical model, the
evolution depends on a fixed parameter $\lambda\in(0, +\infty)$ and
is as follows: at each moment, an infected site becomes healthy at rate
$1$ while a healthy site becomes infected at a rate equal to $\lambda$
times the number of its infected neighbors. For the contact process in
random environment, the single infection parameter~$\lambda$ is
replaced by a collection $(\lambda_e)_{e\in\Ed}$ of random variables
indexed by the set~$\Ed$ of edges of the lattice $\Zd$: the random
variable $\lambda_e$ gives the infection rate between the extremities
of edge $e$, while each site becomes healthy at rate~$1$. We assume
that the law of $(\lambda_e)_{e\in\Ed}$ is stationary and ergodic. From
the application point of view, allowing a random infection rate can be
more realistic in modelizing real epidemics; note that in his
book \cite{MR940469}, Durrett already underlined the inadequacies of
the classical contact process in the modelization of an infection among
a racoon rabbits population, and proposed the contact process in random
environment as an alternative.

Our main result is the following: if we assume that the minimal value
taken by the $(\lambda_e)_{e\in\Ed}$ is above\vadjust{\goodbreak} $\lambda_c(\Zd)$
(the\vspace*{1pt} critical
parameter for the ordinary contact process on $\Zd$), then there exists
a norm $\mu$ on $\Rd$ such that for almost every environment $\lambda
=(\lambda_e)_{e\in\Ed}$, the set $H_t$ of points already infected
before time $t$ satisfies
\[
\Pbarre_{\lambda}\bigl(\exists T>0, t\ge T\Longrightarrow
(1-\epsilon
)tA_{\mu}\subset\tilde{H_t}+ \subset(1+\epsilon)tA_{\mu}\bigr)=1,
\]
where $\tilde{H_t}=H_t+[0,1]^d$, $A_{\mu}$ is the unit ball for $\mu$
and $\Pbarre_{\lambda}$ is the law of the contact process in the
environment $\lambda$, conditioned to survive. The growth of the
contact process in random environment conditioned to survive is thus
asymptotically linear in time, and governed by a shape theorem, as in
the case of the classical contact process on $\Zd$.

Until now, most of the work devoted to the study of the contact process
in random environment focuses on determining conditions for its survival
Liggett~\cite{MR1159569}, Andjel \cite{MR1203176}, Newman and
Volchan \cite{MR1387642} or its extinction Klein~\cite{MR1303643}. They
also mainly deal with the case of dimension $d=1$. Concerning the speed
of the growth when $d=1$, Bramson, Durrett and Schonmann
\cite{MR1112403} show that a random environment can give birth to a
sublinear growth. On the contrary, they conjecture that the growth
should be of linear order for $d\ge2$ as soon as the survival is
possible and that an asymptotic shape result should hold.

For the classical contact process, the proof of the shape result mainly
falls in two parts:
\begin{itemize}
\item The result is first proved for large values of the infection rate
$\lambda$ by Durrett and Griffeath \cite{MR656515} in 1982. They first
obtain, for large $\lambda$, estimates essentially implying that the
growth is of linear order, and then they get the shape result with
superconvolutive techniques.
\item Later, Bezuidenhout and Grimmett \cite{MR1071804} show that a
supercritical contact process conditioned to survive, when seen on a
large scale, stochastically dominates a two-dimensional supercritical
oriented percolation; this guarantees at least linear growth of the
contact process. They also indicate how their construction could be
used to obtain a shape theorem. This last step essentially consists of
proving that the estimates needed in \cite{MR656515} hold for the whole
supercritical regime, and is done by Durrett \cite{MR1117232} in 1989.
\end{itemize}
Similarly, in the case of a random environment, proving a shape theorem
can also fall into two different parts. The first one, and undoubtedly
the hardest one, would be to prove that the growth is of linear order,
as soon as survival is possible; this corresponds to the Bezuidenhout
and Grimmett result in random environment. The second one, which we
tackle here, is to prove a~shape theorem under conditions assuring that
the growth is of linear order; this is the random environment analogous
to the Durrett and Griffeath work. We thus chose to put conditions on
the random environment that allow it to obtain, with classical
techniques, estimates similar to the ones needed in~\cite{MR656515} and
to focus on the proof of the shape result, which already presents
serious additional difficulties when compared to the proof in the
classical case.

The history of shape theorems for random growth models begins in 1961
with Eden \cite{MR0136460} asking for a shape theorem for a tumor
growth model. Richardson \cite{MR0329079} then proves in 1973 a shape
result for a class of models, including Eden model, by using the
technique of subadditive processes initiated in 1965 by Hammersley and
Welsh \cite{MR0198576} for first-passage percolation. From then,
asymptotic shape results for random growth models are usually proved
with the theory of subadditive processes, and, more precisely, with
Kingman's subadditive ergodic theorem \cite{MR0356192} and its
extensions. The most famous example is the shape result for first
passage-percolation on $\Zd$ (see also different variations of this
model Boivin \cite{MR1074741}, Garet and Marchand~\cite{MR2085613},
Vahidi-Asl and Wierman \cite{MR1166620}, Howard and Newmann
\cite{MR1452554}, Howard \cite{MR2023652}, Deijfen
\cite{MR1970474}).

The random growth models can be classified in two families. The first
and most studied one is composed of the permanent models, in which the
occupied set at time $t$ is nondecreasing and extinction is impossible.
First of all are, of course, Richardson models \cite{MR0329079}. More
recently, we can cite the frog model,
introduced in its continuous time version by Bramson and Durrett, and
for which Ram{\'{\i}}rez and Sidoravicius \cite{MR2060478} obtained a
shape theorem, and also the discrete time version, first studied by
Telcs and Wormald \cite{MR1742145} and for which the shape theorem has
been obtained by Alves et al. \cite{MR1910638,MR1893139}. We can also
cite the branching random walks by Comets and Popov \cite{MR2303944}.
In these models, the main part of the work is to prove that the growth
is of linear order, and the whole convergence result is then obtained
by subadditivity.

The second family contains nonpermanent models, in which extinction is
possible. In this case, we rather look for a shape result under
conditioning by the survival. Hammersley \cite{MR0370721} himself, from
the beginning of the subadditive theory, underlined the difficulties
raised by the possibility of extinction. Indeed, if we want to prove
that the hitting times $(t(x))_{x\in\Zd}$ are such that $t(nx)/n$
converges, Kingman's theory requires subadditivity, stationarity and
integrability properties for the collection $t(x)$. Of course, as soon
as extinction is possible, the hitting times can be infinite. Moreover,
conditioning on the survival can break independence, stationarity and
even subadditivity properties. The theory of superconvolutive
distributions was developed to treat cases where either the
subadditivity or the stationarity property lacks; see the lemma
proposed by Kesten in the discussion of Kingman's paper
\cite{MR0356192}, and slightly improved by Hammersley
\cite{MR0370721}, page
674. Note that recently, Kesten and Sidoravicius \cite{MR2247840} use
the same kind of techniques as an ingredient to prove a shape theorem
for a model of the spread of an infection.\looseness=-1

Following Bramson and Griffeath \cite{MR606980,MR578279}, it is on
these ``superconvolutive'' techniques that Durrett and Griffeath
\cite{MR656515} rely to prove the shape result for the classical contact
process on $\Zd$; see also Durrett \cite{MR940469}, that corrects or
clarifies some points of \cite{MR656515}. However, as noticed by
Liggett in the Introduction of \cite{MR806224},
superconvolutive techniques require some kind of independence of the
increments of the process that can limit its application. It is
particularly the case in a random environment setting; for the hitting
times, we have a subadditive property of type
\[
t^\lambda\bigl((n+p)x\bigr)\le t^\lambda(nx)+\tilde{t}^{nx.\lambda}(px)+r(n,p,x).
\]
Here, the exponent gives the environment, $\tilde{t}^{nx.\lambda}(px)$
the same law as the hitting time of $px$ but in the translated
environment $nx.\lambda$, and $r(n,p,x)$ are to be thought of as a
small error term. Following the superconvolutive road would require
that $t^\lambda(nx)$ and $\tilde{t}^{nx.\lambda}(px)$ are independent
and that $\tilde{t}^{nx.\lambda}(px)$ has the same law as
${t}^{\lambda}(px)$.
Now, if we work with a given (quenched) environment, we lose all the
spatial stationarity properties; $\tilde{t}^{nx.\lambda}(px)$ has no
reason to have the same law as ${t}^{\lambda}(px)$. But if we work
under the annealed probability, we lose the markovianity of the contact
process and the independence properties it offers. We thus cannot use,
at least directly, the superconvolutive techniques.

Liggett's extension \cite{MR806224} of the subadditive ergodic theorem
provides an alternate approach when independence properties fail.
However, it does not give the possibility to deal with an error term.
Some works in the same decade (see, e.g., Derriennic
\cite{MR704553}, Derriennic and Hachem \cite{MR939537}
and Sch\"{u}rger \cite{MR833959,MR1127716}) propose almost subadditive ergodic
theorems that do not require independence, but stationarity assumptions
on the extra term are too strong to be used here.
Thus we establish, with techniques inspired from Liggett, a general
subadditive ergodic theorem allowing an error term that matches our situation.

In fact, we do not apply this almost subadditive ergodic theorem
directly to the collection of hitting times $t(x)$, but we rather
introduce the quantity~$\sigma(x)$, that can be seen as a regeneration
time, and that represents a~time when site $x$ is occupied and has
infinitely many descendants. This~$\sigma$ has stationarity and almost
subadditive properties that~$t$ lacks and thus fits the requirements of
our almost subadditive ergodic theorem. Finally, by showing that the
gap between $t$ and $\sigma$ is not too large, we transpose to $t$ the
shape result obtained for $\sigma$.

\section{Model and results}
\label{construction}

\subsection{Environment}

In the following, we denote by \mbox{$\|\cdot\|_1$} and
\mbox{$\|\cdot\|_{\infty}$} the norms on $\Rd$, respectively, defined
by $\|x\|_1=\sum_{i=1}^d |x_i|$ and $\| x\|_{\infty}=\break\max_{1\le i\le d}
|x_i|$. The notation \mbox{$\|\cdot\|$} will be used for an unspecified norm.

We fix $\lambda_c(\Zd)<\lambda_{\min}\le\lambda_{\max}<+\infty
$, where
$\lambda_c(\Zd)$ stands for the critical parameter for the classical
contact process in $\Zd$. Additionally,\vspace*{1pt}\vadjust{\goodbreak} we restrict our study to random
environments $\lambda=(\lambda_e)_{e \in\Ed}$ taking their value in $
\Lambda=[\lambda_{\min},\lambda_{\max}]^{\Ed}$. An environment is thus
a collection $\lambda=(\lambda_e)_{e \in\Ed} \in\Lambda$.

Let $\lambda\in\Lambda$ be fixed. The contact process $(\xi
_t)_{t\ge
0}$ in environment $\lambda$ is a~homogeneous Markov process taking its
values in the set $\mathcal{P}(\Zd)$ of subsets of~$\Zd$. For $z \in
\Zd
$ we also use the random variable $\xi_t(z)=\one_{\{z \in\xi_t\}}$.
If $\xi_t(z)=1$, we say that $z$ is occupied or infected, while if
$\xi
_t(z)=0$, we say that $z$ is empty or healthy. The evolution of the
process is as follows:
\begin{itemize}
\item an occupied site becomes empty at rate $1$,
\item an empty site $z$ becomes occupied at rate
$ \sum_{\|z-z'\|_1=1} \xi_t(z')\lambda_{\{z,z'\}}$,
\end{itemize}
each of these evolutions being independent from the others. In the
following, we denote by $\D$ the set of c\`{a}dl\`{a}g functions from
$\R_{+}$ to $\mathcal{P}(\Zd)$; it is the set of trajectories for
Markov processes with state space $\mathcal{P}(\Zd)$.

To define the contact process in environment $\lambda\in\Lambda$, we
use the Harris construction \cite{MR0488377}. It allows us to couple
contact processes starting from distinct initial configurations by
building them from a single collection of Poisson measures on~$\R_+$.

\subsection{Construction of the Poisson measures}

We endow $\R_+$ with the Borel $\sigma$-algebra $\mathcal B(\R_+)$, and
we denote by $M$ the set of locally finite counting measures $m=\sum
_{i=0}^{+\infty} \delta_{t_i}$. We endow this set with the $\sigma
$-algebra $\mathcal M$ generated by the maps $m\mapsto m(B)$, where $B$
describes the set of Borel sets in $\R_+$.

We then define the measurable space $(\Omega, \mathcal F)$ by setting
\[
\Omega=M^{\Ed}\times M^{\Zd} \quad\mbox{and}\quad \mathcal F=\mathcal
{M}^{\otimes\Ed} \otimes\mathcal{M}^{\otimes\Zd}.
\]
On this space, we consider the family $(\P_{\lambda})_{\lambda\in
\Lambda
}$ of probability measures defined as follows:
for every $\lambda=(\lambda_e)_{e \in\Ed} \in\Lambda$,
\[
\P_{\lambda}=\biggl(\bigotimes_{e \in\Ed} \mathcal{P}_{\lambda
_{e}}\biggr) \otimes\mathcal{P}_1^{\otimes\Zd},
\]
where, for every $\lambda\in\R_+$, $\mathcal{P}_{\lambda}$ is the law
of a punctual Poisson process on $\R_+$ with intensity $\lambda$. If
$\lambda\in\R_+$, we write $\P_\lambda$ (rather than $\P
_{(\lambda
)_{e \in\Ed}}$) for the law in deterministic environment with constant
infection rate $\lambda$.

For every $t\ge0$, we denote by $\mathcal{F}_t$ the $\sigma$-algebra
generated by the maps $\omega\mapsto\omega_e(B)$ and $\omega\mapsto
\omega_z(B)$, where $e$ ranges over all edges in $\Ed$, $z$ ranges over
all points in $\Zd$ and $B$ ranges over the set of Borel sets in $[0,t]$.

\subsection{Graphical construction of the contact process}

This construction is exposed in all details in Harris \cite{MR0488377};
we just give here an informal description. Let $\omega=((\omega_e)_{e
\in\Ed}, (\omega_z)_{z \in\Zd}) \in\Omega$. Above each site $z
\in
\Zd$, we draw a time line~$\R_+$, and we put a cross at the times given
by $\omega_z$. Above each edge $e \in\Ed$, we draw at the times given
by $\omega_e$ a horizontal segment between the extremities of the edge.\vadjust{\goodbreak}

An open path follows the time lines above sites (but crossing crosses
is forbidden) and uses horizontal segments to jump from a time line to
a neighboring time line; in this description, the evolution of the
contact process looks like a percolation process, oriented in time but
not in space.
For $x,y \in\Zd$ and $t \ge0$, we say that $\xi_t^x(y)=1$ if and only
if there exists an open path from $(x,0)$ to $(y,t)$, then we define
%
%
\begin{eqnarray} \label{additivite}
\xi_t^x & = & \{y \in\Zd\dvtx\xi_t^x(y)=1\}, \nonumber\\[-8pt]\\[-8pt]
\forall A \in\mathcal P(\Zd)\qquad \xi_t^A & = & \bigcup_{x \in A}
\xi
_t^x.\nonumber
\end{eqnarray}
For instance, we obtain
$(A \subset B) \Rightarrow(\forall t \ge0$ $\xi_t^A \subset\xi_t^B)$.

When $\lambda\in\R_{+}^*$, Harris shows that under $\P_{\lambda}$, the
process $(\xi^A_t)_{t \ge0}$ is the contact process with infection
rate $\lambda$, starting from initial configuration $A$. The proof can
readily be extended to a nonconstant $\lambda\in\Lambda$, which
allows us to define the contact process in environment $\lambda$
starting from initial configuration $A$. This is a Feller process, and
thus it benefits from the strong Markov property.

\subsection{Time translations}

For $t \ge0$, we define the translation operator $\theta_t$ on a
locally finite counting measure $m=\sum_{i=1}^{+\infty} \delta_{t_i}$
on $\R_+$ by setting
\[
\theta_t m=\sum_{i=1}^{+\infty} \one_{\{t_i\ge t\}}\delta_{t_i-t}.
\]
The translation $\theta_t$ induces an operator on $\Omega$, still
denoted by $\theta_t$; for every $\omega\in\Omega$, we set
\[
\theta_t \omega=((\theta_t \omega_e)_{e \in\Ed}, (\theta_t
\omega
_z)_{z \in\Zd}).
\]
The Poisson point process being translation invariant, every
probability measure $\P_{\lambda}$ is stationary under $\theta_t$. The
semigroup property of the contact process here has a stronger
trajectorial version; for every $A \subset\Zd$, for every $s,t \ge0$,
for every $\omega\in\Omega$, we have
%
%
\begin{equation}
\label{semigroupe}
\xi_{t+s}^A(\omega) = \xi_{s}^{\xi_t^A(\omega)}(\theta_t\omega
)=\xi
_{s}^{\centerdot}(\theta_t\omega)\circ\xi_t^A(\omega),
\end{equation}
that can also be written in the classical markovian way
\[
\forall B\in\mathcal{B}(\mathcal{D})\qquad \P\bigl((\xi_{t+s}^A)_{s \ge
0} \in
B| \mathcal F_t\bigr) = \P\bigl((\xi_{s}^{\centerdot})_{s \ge0} \in B\bigr) \circ
\xi^A_t.
\]
We can write in the same way the strong Markov property: if $T$ is an
$(\mathcal{F}_t)_{t\ge0}$ stopping time, then, on the event $\{
T<+\infty\}$,
\begin{eqnarray*}
\xi_{T+s}^A(\omega)& = & \xi_{s}^{\xi_T^A(\omega)}(\theta_T\omega
),\\
\forall B\in\mathcal{B}(\mathcal{D})\qquad \P\bigl((\xi_{T+s}^A)_{s
\ge0} \in B| \mathcal F_T\bigr) & = & \P\bigl((\xi_{s}^{\centerdot})_{s \ge0}
\in B\bigr) \circ\xi^A_T.
\end{eqnarray*}
We recall that $\mathcal{F}_T$ is defined by
\[
\mathcal{F}_T=\bigl\{B\in\mathcal{F}\dvtx\forall t\ge0 B\cap\{
T\le
t\}\in\mathcal{F}_t\bigr\}.
\]

\subsection{Spatial translations}

The group $\Zd$ can act on the process and on the environment. The
action on the process changes the observer's point of view of the process.
For $x \in\Zd$, we define the translation operator $T_x$ by
\[
\forall\omega\in\Omega\qquad T_x \omega=(( \omega_{x+e})_{e \in
\Ed},
( \omega_{x+z})_{z \in\Zd}),
\]
where $x+e$ the edge $e$ translated by vector $x$.

Besides, we can consider the translated environment $\T{x}{\lambda}$
defined by $(\T{x}{\lambda})_e=\lambda_{x+e}$.
These actions are dual in the sense that for every $\lambda\in\Lambda
$, for every $x \in\Zd$,
%
%
\begin{equation}
\label{translationspatiale}
\forall A\in\mathcal{F} \qquad\P_{\lambda}(T_x \omega\in A) =
\P_{\T
{x}{\lambda}}(\omega\in A).
\end{equation}
Consequently, the law of $\xi^x$ under $\P_\lambda$ coincides with the
law of $\xi^0$ under~$\P_{x.\lambda}$.

\subsection{Essential hitting times and associated translations}

For a set \mbox{$A \subset\Zd$}, we define the life time $\tau^A$ of the
process starting from $A$ by
\[
\tau^A=\inf\{t\ge0\dvtx\xi_t^A=\varnothing\}.
\]
For $A \subset\Zd$ and $x \in\Zd$, we also define the first infection
time $t^A(x)$ of site~$x$ from set $A$ by
\[
t^A(x)=\inf\{t\ge0\dvtx x \in\xi_t^A\}.
\]
If $y\in\Zd$, we write $t^y(x)$ instead of $t^{\{y\}}(x)$. Similarly,
we simply write $t(x)$ for $t^0(x)$.

We now introduce the essential hitting time $\sigma(x)$: it is a time
where the site $x$ is infected from the origin $0$ and also has an
infinite life time. This essential hitting time is defined through a
family of stopping times as follows. We set $u_0(x)=v_0(x)=0$ and we
define recursively two increasing sequences of stopping times
$(u_n(x))_{n \ge0}$ and $(v_n(x))_{n \ge0}$ with
$u_0(x)=v_0(x)\le u_1(x)\le v_1(x)\le u_2(x)\cdots$ as follows:
\begin{itemize}
\item Assume that $v_k(x)$ is defined. We set $u_{k+1}(x) =\inf\{t\ge
v_k(x)\dvtx x \in\xi^0_t \}$.
If $v_k(x)<+\infty$, then $u_{k+1}(x)$ is the first time after $v_k(x)$
where site $x$ is once again infected; otherwise, $u_{k+1}(x)=+\infty$.
\item Assume that $u_k(x)$ is defined, with $k \ge1$. We set
$v_k(x)=u_k(x)+\tau^x\circ\theta_{u_k(x)}$.
If $u_k(x)<+\infty$, the time $\tau^x\circ\theta_{u_k(x)}$ is the life
time of the contact process starting from $x$ at time $u_k(x)$;
otherwise, $v_k(x)=+\infty$.
\end{itemize}
We then set
%
%
\begin{equation}
\label{definitiondeK}
K(x)=\min\{n\ge0\dvtx v_{n}(x)=+\infty\mbox{ or } u_{n+1}(x)=+\infty
\}.
\end{equation}
This quantity represents the number of steps before the success of this
process; either we stop because we have just found an infinite
$v_n(x)$, which corresponds to a time $u_n(x)$ when $x$ is occupied and
has infinite progeny, or we stop because we have just found an infinite
$u_{n+1}(x)$, which says that after $v_n(x)$, site $x$ is never
infected anymore.

We then set $\sigma(x)=u_{K(x)}$, and call it the essential hitting
time of $x$. It is, of course, larger than the hitting time $t(x)$ and
can been seen as a regeneration time.

Note however that $\sigma(x)$ is not necessary the first time when $x$
is occupied and has infinite progeny. For instance, such an event can
occur between~$u_1(x)$ and~$v_1(x)$, being ignored by the recursive
construction.

We will see that $K(x)$ is almost surely finite, so $\sigma(x)$ is
well defined.
At the same time, we define the operator $\tilde\theta_x$ on $\Omega
$ by
\[
\tilde\theta_x =
\cases{
T_{x} \circ\theta_{\sigma(x)}, &\quad if $\sigma(x)<+\infty$,\cr
T_x, &\quad otherwise,}
\]
or, more explicitly,
\[
(\tilde\theta_x)(\omega) =
\cases{
T_{x} \bigl(\theta_{\sigma(x)(\omega)} \omega\bigr), &\quad if
$\sigma(x)(\omega)<+\infty$,\cr
T_x (\omega), &\quad otherwise.}
\]
We will mainly deal with the essential hitting time $\sigma(x)$ that
enjoys, unlike~$t(x)$,
some good invariance properties in the survival-conditioned
environment. We will also control the difference between $\sigma(x)$
and $t(x)$, which will allow us to transpose to $t(x)$ the results
obtained for $\sigma(x)$.

\subsection{Contact process in the survival-conditioned environment}

We now have to introduce the random environment. In the following, we
fix a probability measure $\nu$ on the sets of environments $\Lambda
=[\lambda_{\min},\lambda_{\max}]^{\Ed}$. We assume that $\nu$ is
stationary and, denoting by $\Erg(\nu)$ the set of $x\in\Zd \setminus\{
0\}$ such that the translation by $x$ is ergodic for $\nu$, then the
cone generated by $\Erg(\nu)$ is dense in $\Rd$. This condition is
obviously fulfilled if $\Erg(\nu)=\Zd\setminus\{ 0\}$. This perhaps odd
condition allows us to consider some natural models where the
ergodicity assumption is not satisfied in some directions, for example,
along the coordinate vectors. This setting naturally contains the case
of an i.i.d. random environment and the case\vspace*{1pt} of a
deterministic environment $\lambda>\lambda_c(\Zd)$; we simply take for
$\nu$ the Dirac measure $(\delta_{\lambda})^{\otimes\Ed}$.

For $\lambda\in\Lambda$, we define the probability measure ${\Pbarre
}_\lambda$
on $(\Omega, \mathcal F)$ by
\[
\forall E\in\mathcal{F}\qquad {\Pbarre}_\lambda(E)=\P_\lambda
(E|\tau
^0=+\infty).
\]
It is thus the law of the family of Poisson point processes,
conditioned to the survival of the contact process starting from $0$.
On the same space $(\Omega, \mathcal F)$, we define the corresponding
annealed probability $\Pbarre$ by setting
\[
\forall E\in\mathcal{F} \qquad{\Pbarre}(E)=\int_\Lambda{\Pbarre
}_\lambda
(E) \,d\nu(\lambda).
\]
In other words, the environment $\lambda=(\lambda_e)_{e \in\Ed}$ where
the contact process lives is a random variable with law $\nu$, and it
is under the probability measure $\Pbarre$ that we seek the asymptotic
shape theorem.\vadjust{\goodbreak}

It could seem more natural to work with the following probability measure:
\[
\forall E\in\mathcal{F}\qquad \hat{\P}(E)=\P(E|\tau^0=+\infty)=
\frac{\int\Pbarre_{\lambda}(E)\P_\lambda(\tau^0=+\infty) \,d\nu
(\lambda
)}{\int\P_\lambda(\tau^0=+\infty) \,d\nu(\lambda)}.
\]
It appears that our proofs do not work with this probability measure.
However, our restrictions on the set $\Lambda$ of possible environments
ensure that~$\Pbarre$ and~$\hat{\P}$ are equivalent; the $
\Pbarre$-a.s.
asymptotic shape theorem is thus also a~$\hat{\P}$-a.s. asymptotic
shape theorem.

\subsection{Organization of the paper and results}

In Section \ref{sigma}, we establish the invariance and ergodicity
properties. In particular, we prove the following theorem.
%
%
\begin{theorem}
\label{systemeergodique}
For every $x\in\Erg(\nu)$, the measure-preserving dynamical system
$(\Omega,\mathcal{F},\Pbarre,\tilde{\theta}_x)$ is ergodic.
\end{theorem}

In Section \ref{controlediff}, we study the integrability properties of
the family $(\sigma(x))_{x\in\Zd}$; we also control the discrepancy
between $\sigma(x)$ and $t(x)$ and the lack of subadditivity
of $\sigma$.
%
%
\begin{theorem}
\label{presquesousadditif}
There exist $A_{\rref{epresquesousadditif}},B_{\rref
{epresquesousadditif}}>0$ such that
for any $\lambda\in\Lambda$, for any $x$, \mbox{$y\in\Zd$},
%
%
\begin{equation}\label{epresquesousadditif}\quad
\forall t>0 \qquad\Pbarre_\lambda\bigl(\sigma(x+y)-\bigl(\sigma(x)+\sigma(y)\circ
\tilde{\theta}_x\bigr)\ge t\bigr)\le A_{\theequation}\exp\bigl(-B_{\theequation
}\sqrt{t}\bigr).
\end{equation}
\end{theorem}

At first sight, one could think that $\sigma(x+y)\le\sigma(x)+\sigma
(y)\circ\tilde{\theta}_x$ always holds, but this is not the case
because $\sigma(x+y)$ is not necessary the first time when $x+y$ is
occupied and has infinite progeny.

However, the theorem says that the lack of subadditivity of $\sigma$ is
really small; in particular, it does not depend on the considered points.
Then, in the same spirit as Kingman \cite{MR0438477} and Liggett
\cite{MR806224}, we prove in Section \ref{forme} that for every $x\in
\Zd$,
the ratio $\frac{\sigma(nx)}n$ converges $\Pbarre$-a.s. to a real
number $\mu(x)$.
The functional $x\mapsto\mu(x)$ can be extended into a norm on $\Rd$,
which will characterize the asymptotic shape. In the following, $A_{\mu
}$ will denote the unit ball for $\mu$. We define the sets
\begin{eqnarray*}
H_t & = & \{x\in\Zd\dvtx t(x)\le t\},\\
G_t & = & \{x\in\Zd\dvtx\sigma(x)\le t\},\\
K'_t & = & \{x\in\Zd\dvtx\forall s\ge t\mbox{ }\xi^0_s(x)=\xi^{\Zd
}_s(x)\},
\end{eqnarray*}
and we denote by $\tilde{H}_t,\tilde{G}_t,\tilde{K}'_t$ their
``fattened'' versions
\[
\tilde{H}_t=H_t+[0,1]^d, \qquad\tilde{G}_t=G_t+[0,1]^d
\quad\mbox{and}\quad
\tilde
{K}'_t=K'_t+[0,1]^d.
\]
We can now state the asymptotic shape result.
%
%
\begin{theorem}[(Asymptotic shape theorem)]
\label{thFA}
For every $\epsilon>0$, $\Pbarre$-a.s., for every $t$ large enough,
%
%
\begin{equation}
\label{leqdeforme}
(1-\epsilon)A_{\mu}\subset\frac{\tilde K'_t\cap\tilde G_t}t\subset
\frac{\tilde G_t}t\subset\frac{\tilde H_t}t\subset(1+\epsilon
)A_{\mu}.
\end{equation}
\end{theorem}

The set $K'_t\cap G_t$ is the coupled zone of the process. Usually, the
asymptotic shape result for the coupled zone is rather expressed in
terms of $K_t\cap H_t$,
where
\[
K_t=\{x\in\Zd\dvtx\xi^0_t(x)=\xi^{\Zd}_t(x)\}.
\]
Our result also gives the shape theorem for $K_t\cap H_t$, because
$K'_t\cap G_t\subset K_t\cap H_t\subset H_t$.

Let us note that the shape result can also be formulated in the
following ``quenched'' terms: for $\nu$-a.e. environment, we know that
on the event
``the contact process survives,'' its growth is governed by
(\ref{leqdeforme}) for $t$ large enough.
We can also give a complete convergence result.
%
%
\begin{theorem}[(Complete convergence theorem)]
\label{thCC}
For every $\lambda\in\Lambda$, the contact process in environment
$\Lambda$ admits an upper invariant measure $m_\lambda$ defined by
\[
\forall A\subset\Zd, \Card{A}<+\infty\qquad m_\lambda(\omega
\supset
A)=\lim_{t\to+\infty}\P_{\lambda}(\xi^{\Zd}_t\supset A).
\]
Then, for every finite set $A\subset\Zd$ and for $\nu$-a.e.
environment $\lambda$, one has
\[
\P^A_{\lambda,t} \Longrightarrow\P_{\lambda}(\tau^A<+\infty
)\delta
_{\varnothing}+\P_{\lambda}(\tau^A=\infty)m_{\lambda},
\]
where $\P^A_{\lambda,t}$ is the law of $\xi^A_t$ under $\P_{\lambda}$
and $\Longrightarrow$ stands for the convergence in law.
\end{theorem}

The proof of this result does not require any new idea, and we just
give a hint at the end of Section \ref{restartestimees}.

As explained in the \hyperref[intro]{Introduction}, in order to prove
the asymptotic
shape theorem, we need some estimates analogous to the ones needed in
the proof by Durrett and Griffeath in the classical case. We set
\[
B_r^x=\{y \in\Zd\dvtx\|y-x\|_{\infty} \le r\},
\]
and we write $B_r$ instead of $B_r^0$.
%
%
\begin{prop}
\label{propuniforme}
There exist $A,B,M,c,\rho>0$ such that for every
$\lambda\in\Lambda$, for every $y \in\Zd$, for every $ t\ge0$,
%
%
\begin{eqnarray}
\label{uniftau}
\P_\lambda(\tau^0=+\infty) & \ge& \rho,
\\
\label{richard}
\P_\lambda(H^0_t \not\subset B_{Mt} ) & \le& A\exp(-Bt),
\\
\label{grosamasfinis}
\P_\lambda( t<\tau^0<+\infty) &\le& A\exp(-Bt),
\\
\label{retouche}
\P_{\lambda}\biggl( t^0(y)\ge\frac{\|y\|}c+t, \tau^0=+\infty\biggr)
& \le& A\exp(-Bt),
\\
\label{petitsouscouple}
\P_{\lambda}(0\notin K'_t, \tau^0=+\infty) &\le&A\exp(-B t).
\end{eqnarray}
\end{prop}

All these estimates are already available for the classical contact
process in the supercritical regime. For large $\lambda$, they are
established by Durrett and Griffeath~\cite{MR656515}, and the extension
to the entire supercritical regime is made possible thanks to
Bezuidenhout and Grimmett's work \cite{MR1071804}. For the crucial
estimate (\ref{grosamasfinis}), one can find the detailed proof in
Durrett \cite{MR1117232} or in Liggett~\cite{MR1717346}. The need for
these estimates explains our restrictions on the possible range of the
random environment.

We chose to focus on the stationarity and subadditivity properties of
the essential hitting time $\sigma$ and on the proof of the shape
result. We thus admit in Sections \ref{sigma}, \ref{controlediff}
and \ref{forme} the uniform controls given by Proposition \ref
{propuniforme}, whose proof (via restart arguments)
is postponed to Section \ref{restartestimees}. That section is
totally independent of the rest of the paper. Finally, in
the \hyperref[app]{Appendix}, we prove a general (almost) subadditive
ergodic theorem. As we think it could also be useful in other
situations, we present it in a more general form than what is needed
for our aim.

\section{\texorpdfstring{Properties of $\tilde{\theta}_x$}{Properties of theta x}}
\label{sigma}

\subsection{First properties}

We first check that $K(x)$ is almost surely finite and even has a
subgeometrical tail.
%
%
\begin{lemme}
\label{Kgeom}
$ \forall A\subset\Zd,\forall x\in\Zd,
\forall
\lambda\in\Lambda,\forall n\in\N$ $\P_\lambda(K(x)>n)\le
(1-\rho)^n$.
\end{lemme}
\begin{pf}
Remember that $\rho$ is given in (\ref{uniftau}). Let $\lambda\in
\Lambda$ and $n \in\N$. The strong Markov property applied at time
$u_{n+1}(x)$ ensures that
\begin{eqnarray*}
\P_\lambda\bigl(K(x)>n+1\bigr)
& = & \P_\lambda\bigl(u_{n+2}(x)<+\infty\bigr) \\
& \le& \P_\lambda\bigl(u_{n+1}(x)<+\infty,v_{n+1}(x)<+\infty\bigr) \\
& \le& \P_\lambda\bigl(u_{n+1}(x)<+\infty,\tau^x \circ\theta
_{u_{n+1}(x)}<+\infty\bigr) \\
& \le& \P_\lambda\bigl(u_{n+1}(x)<+\infty\bigr)\P_\lambda(\tau^x <+\infty
)\\
& \le& \P_\lambda\bigl(u_{n+1}(x)<+\infty\bigr) (1-\rho)\\
&=&\P_\lambda
\bigl(K(x)>n\bigr)(1-\rho),
\end{eqnarray*}
which proves the lemma.
\end{pf}
%
%
\begin{lemme}
Let $\lambda\in\Lambda$. $\P_\lambda$-a.s., for every $x\in\Zd$,
%
%
\begin{eqnarray}
\label{bleue}
&&\bigl(K(x)=k\bigr)\quad\mbox{and}\nonumber\\[-8pt]\\[-8pt]
&&(\tau^0=+\infty) \quad
\Longleftrightarrow\quad
\bigl(u_k(x)<+\infty\mbox{ and } v_k(x)=+\infty\bigr).\nonumber
\end{eqnarray}
\end{lemme}
\begin{pf}
Let $\lambda\in\Lambda$. By Lemma \ref{Kgeom}, the number $K(x)$ is
$\P_{\lambda}$-a.s. finite. Let $k \in\N$; the strong Markov property
applied at time $v_k(x)$ ensures that
\begin{eqnarray*}
& & \P_{\lambda}\bigl(\tau^0=+\infty, v_k(x)<+\infty,
u_{k+1}(x)=+\infty|
\mathcal F_{v_k(x)}\bigr)\\
&&\qquad = \one_{\{v_k(x)<+\infty\}} \P_\lambda\bigl(\tau^{\centerdot
}=+\infty,
t^{\centerdot}(x)=+\infty\bigr) \circ\xi^0_{v_k(x)}.
\end{eqnarray*}
Consider now a finite nonempty set $B\subset\Zd$. With (\ref
{retouche}), we get
\begin{eqnarray*}
\P_\lambda\bigl( \tau^B=+\infty, t^B(x)=+\infty\bigr)
& \le& \sum_{y \in B}\P_\lambda\bigl(\tau^y=+\infty, t^y(x)=+\infty\bigr)
\\
& \le& \sum_{y \in B}\P_{\T{y}{\lambda}}\bigl(\tau^0=+\infty,
t^0(x-y)=+\infty\bigr)=0.
\end{eqnarray*}
This gives the direct implication. The reverse one comes from (\ref
{semigroupe}).
\end{pf}

Our construction of $\sigma(x)$ is very similar to the restart process
exposed in Durrett and Griffeath \cite{MR656515}. The essential
difference is that in that paper, the aim is to find, close to $x$, a
point that survives while we require here the point to be exactly at
$x$. Thus, we will be able to describe precisely the law of the contact
process starting from $x$ at time $\sigma(x)$, and construct
transformations under which $\Pbarre$ is invariant.
%
%
\begin{lemme}
\label{magic}
Let $x \in\Zd\setminus\{0\}$, $A$ in the $\sigma$-algebra generated
by $\sigma(x)$ and $B\in\mathcal F$. Then
\[
\forall\lambda\in\Lambda\qquad\Pbarre_\lambda\bigl(A \cap(\tilde
{\theta
}_x)^{-1}(B)\bigr)=\Pbarre_\lambda(A) \Pbarre_{\T{x}{\lambda}}(B).
\]
\end{lemme}
\begin{pf}
We just have to check that for any $k\in\N^*$, one has
\[
\Pbarre_\lambda\bigl(A \cap(\tilde{\theta}_x)^{-1}(B) \cap\{K(x)=k\}
\bigr)=\Pbarre_\lambda\bigl(A\cap\{K(x)=k\}\bigr) \Pbarre_{\T{x}{\lambda}}(B).
\]
Consider a Borel set $A'\subset\R$ such that $A=\{\sigma(x) \in A'\}$.
The essential hitting time $\sigma(x)$ is not a stopping time, but we
can use the stopping times of the construction
%
%
\begin{eqnarray}
&&\P_\lambda\bigl(\{\tau^0=+\infty\} \cap A \cap(\tilde{\theta}_x)^{-1}(B)
\cap\{K(x)=k\}\bigr) \nonumber\\
\label{un}
&&\qquad = \P_\lambda\bigl(\tau^0=+\infty, \sigma(x) \in A', T_x
\circ
\theta_{\sigma(x)} \in B, u_k(x)<+\infty, v_k=+\infty\bigr)
\\
\label{deux}
&&\qquad = \P_\lambda\bigl(u_k(x)<+\infty, u_k(x) \in A', \tau^x\circ
\theta
_{u_k(x)}=+\infty, T_x \circ\theta_{u_k(x)}\in B\bigr)\\
\label{trois}
&&\qquad = \P_\lambda\bigl( u_k(x) \in A', u_{k}(x)<+\infty\bigr) \P_{\lambda
}(\tau^x=+\infty, T_x \in B)\\
\label{quatre}
&&\qquad = \P_\lambda\bigl( u_k(x) \in A', u_{k}(x)<+\infty\bigr)\P_{\T
{x}{\lambda
}}(\{\tau^0=+\infty\} \cap B).
\end{eqnarray}
For (\ref{un}), we use equivalence (\ref{bleue}). For (\ref{deux}), we
notice that for any stopping time~$T$,
%
%
\begin{equation}
\label{violette}
\{T<+\infty, x \in\xi_T^0, \tau^0 \circ T_x \circ\theta
_T=+\infty\}\subset\{\tau^0=+\infty\}.
\end{equation}
Equality (\ref{trois}) follows from the strong Markov property applied
at time $u_{k}(x)$, while (\ref{quatre}) comes from the spatial
translation property (\ref{translationspatiale}). Dividing the identity
by $\P_{\lambda}(\tau^0=+\infty)$, we obtain an identity of the form
\[
\Pbarre_\lambda\bigl(A \cap(\tilde{\theta}_x)^{-1}(B) \cap\{K(x)=k\}
\bigr)=\psi
(x,\lambda,k,A) \Pbarre_{\T{x}{\lambda}}(B),
\]
and the number $\psi(x,\lambda,k,A)$ is identified by taking
$B=\Omega$.
\end{pf}
%
%
\begin{coro}
\label{invariancePbarre}
Let $x,y \in\Zd$ and $\lambda\in\Lambda$. Assume that $x \neq0$.
\begin{itemize}
\item The probability measure $\Pbarre$ is invariant under the
translation $\tilde\theta_x$.
\item Under $\Pbarre_\lambda$, $\sigma(y)\circ\tilde{\theta}_x$ and
$\sigma(x)$ are independent. Moreover, the law of $\sigma(y)\circ
\tilde
{\theta}_x$ under $\Pbarre_\lambda$ is the same as the law of
$\sigma
(y)$ under $\Pbarre_{\T{x}{\lambda}}$.
\item The random variables $(\sigma(x) \circ(\tilde\theta_{x})^j)_{j
\ge0}$ are independent under $\Pbarre_\lambda$.
\end{itemize}
\end{coro}
\begin{pf}
For the first point, we just apply the previous lemma with \mbox{$A=\Omega$},
then we integrate with respect to $\lambda$ and use the stationarity of
$\nu$.

For the second point, let $A',B'$ be two Borel sets in $\R$ and apply
Lemma~\ref{magic} with $A=\{\sigma(x)\in A'\}$ and $B=\{\sigma(y)
\circ
\tilde{\theta}_x\in B'\}$.

Let $n \ge1$ and $A_0,A_1, \ldots, A_n$ be some Borel sets in $\R$.
We have
\begin{eqnarray*}
\hspace*{-4pt}&& \Pbarre_\lambda\bigl(\sigma(x) \in A_0, \sigma(x) \circ\tilde\theta_{x}
\in A_1, \ldots, \sigma(x) \circ(\tilde\theta_{x})^n \in A_n\bigr) \\
\hspace*{-4pt}&&\qquad = \Pbarre_\lambda\bigl(\sigma(x) \in A_0, \bigl(\sigma(x),
\ldots, \sigma(x) \circ(\tilde\theta_{x})^{n-1}\bigr)\circ\tilde\theta
_{x} \in A_1\times\cdots\times A_n\bigr) \\
\hspace*{-4pt}&&\qquad = \Pbarre_\lambda\bigl(\sigma(x) \in A_0\bigr) \Pbarre_{\T{x}{\lambda
}}\bigl(\sigma(x) \in A_1, \sigma(x) \circ\tilde\theta_{x} \in A_2,
\ldots
, \sigma(x) \circ(\tilde\theta_{x})^{n-1} \in A_n\bigr),
\end{eqnarray*}
where the last equality comes from Lemma \ref{magic}.
We recursively obtain
\[
\Pbarre_\lambda\biggl( \bigcap_{0 \le j \le n} \{ \sigma(x) \circ
(\tilde\theta_{x})^j \in A_j\}\biggr)=\prod_{0 \le j \le n} \Pbarre
_{\T
{jx}{\lambda}} \bigl(\sigma(x)\in A_j\bigr),
\]
which ends the proof of the lemma.
\end{pf}

\subsection{Ergodicity}

To prove Theorem \ref{systemeergodique}, it seems natural to estimate
the evolution with $m$ of the dependence between
$A$ and $\tilde{\theta}_x^{-m}(B)$ for some events~$A$ and $B$.
If $m \ge1$, the operator $\tilde{\theta}_x^m$ corresponds to a
spatial translation by vector $mx$ and to a time translation by vector $S_m(x)$
\[
\tilde{\theta}_x^m  = T_{mx} \circ\theta_{S_m(x)}\qquad
\mbox{with } S_m(x)  = \sum_{j=0}^{m-1} \sigma(x) \circ\tilde
{\theta}_x^j.
\]
We begin with a lemma in the same spirit as Lemma \ref{magic}.
%
%
\begin{lemme}
\label{supermagic}
Let $t>0$, $A \in\mathcal F_t$ and $B\in\mathcal F$.

Then, for any $x \in\Zd$, any $\lambda\in\lambda$, any $m \ge1$,
\[
\Pbarre_\lambda\bigl(A \cap\{t\le S_m(x)\} \cap(\tilde{\theta
}_x^m)^{-1}(B)\bigr)=\Pbarre_\lambda\bigl(A\cap\{t\le S_m(x)\}\bigr) \Pbarre_{\T
{mx}{\lambda}}(B).
\]
\end{lemme}
\begin{pf}
Set $\overline{K}_m(x)=(K(x),K(x)\circ\tilde{\theta}_x,\ldots
,K(x)\circ
\tilde{\theta}_x^{m-1})$. It is sufficient to prove that for any
$k=(k_0,\ldots,k_{m-1})\in(\N^*)^m$, one has
\begin{eqnarray*}
& & \Pbarre_\lambda\bigl(A, t\le S_m(x), \tilde{\theta}_x^{-m}(B),
\overline{K}_m(x)=k\bigr) \\
&&\qquad = \Pbarre_\lambda\bigl(A, t\le S_m(x), \overline{K}_m(x)=k\bigr)
\Pbarre
_{\T{mx}{\lambda}}(B).
\end{eqnarray*}

Let $k\in(\N^*)^m$. We set $R_0(x)=0$ and, for $l \le m-2$,
$R_{l+1}(x)=R_l +u_{k_l}(x) \circ\theta_{R_l(x)}$. Thanks to remark
%
%
\begin{figure}

\includegraphics{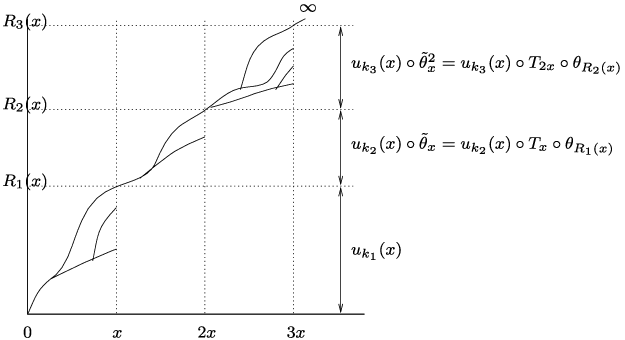}

\caption{An example with $k_1=3$, $k_2=2$ and $k_3=4$.}
\label{uneautrefigure}
\end{figure}
(\ref{violette}), the following events coincide (see Figure
\ref{uneautrefigure}):
\begin{eqnarray*}
&&\{
\tau^0=+\infty,
\overline{K}_m(x)=k
\}\\
&&\qquad=
\bigl\{
u_{k_1}(x)<+\infty,
u_{k_2}(x) \circ T_x\circ\theta_{R_1(x)}<+\infty, \ldots,\\
&&\hspace*{38pt}
u_{k_m}(x) \circ T_{(m-1)x} \circ\theta_{R_{m-1}(x)}<+\infty,
\tau^0 \circ T_{mx}\circ\theta_{R_m(x)}=+\infty
\bigr\}.
\end{eqnarray*}
Moreover, on this event, $S_m(x)=R_m(x)$ holds. Thus
\begin{eqnarray*}
&&\P_\lambda
\bigl(\tau^0=+\infty, A,
t\le S_m(x),
\overline{K}_m(x)=k, \tilde{\theta}_x^{-m}(B)\bigr)
\\
&&\qquad = \P_\lambda
\bigl(A, u_{k_1}(x)<+\infty,
u_{k_2}(x) \circ T_x\circ\theta_{R_1(x)}<+\infty, \ldots,\\
&&\qquad\quad\hspace*{17pt}u_{k_m}(x) \circ T_{(m-1)x} \circ\theta_{R_{m-1}(x)}<+\infty,
t\le R_m(x),\\
&&\qquad\quad\hspace*{49pt}\tau^0 \circ T_{mx}\circ\theta_{R_m(x)}=+\infty,
T_{mx} \circ
\theta_{R_m(x)} \in B\bigr).
\end{eqnarray*}
By construction, $R_m(x)$ is a stopping time and the event
\[
A \cap\{u_{k_1}(x)<+\infty\} \cap\cdots\cap\bigl\{ u_{k_m}(x) \circ
T_{(m-1)x} \circ\theta_{R_{m-1}(x)}<+\infty\bigr\} \cap\{t\le R_m(x)\}
\]
is measurable with respect to $\mathcal F_{R_m(x)}$. Using the strong
Markov property and the spatial translation property (\ref
{translationspatiale}), we get
\begin{eqnarray*}
&&\P_\lambda
\bigl(\tau^0=+\infty, A,
t\le S_m(x),
\overline{K}_m(x)=k, \tilde{\theta}_x^{-m}(B)\bigr)\\
&&\qquad = \P_\lambda
\bigl(A, u_{k_1}(x)<+\infty,
u_{k_2}(x) \circ T_x\circ\theta_{u_{k_1}(x)}<+\infty, \ldots,\\
&&\qquad\quad\hspace*{33.5pt}u_{k_m}(x) \circ T_{(m-1)x} \circ\theta_{R_{m-1}(x)}<+\infty,
t\le R_m(x)\bigr)
\\
&&\qquad\quad{} \times\P_{\T{mx}{\lambda}}(\{\tau=+\infty\} \cap B ).
\end{eqnarray*}
Dividing the identity by $\P_{\lambda}(\tau=+\infty)$, we obtain an
identity of the form
\[
\Pbarre_\lambda\bigl(A, t\le S_m(x), \tilde{\theta}_x^{-m}(B),
\overline{K}_m(x)=k\bigr)=\psi(x,\lambda,k,m,A) \Pbarre_{\T{mx}{\lambda}}(B),
\]
and we identify the value of $\psi(x,\lambda,k,m,A)$ by taking
$B=\Omega$.
\end{pf}

We can now state a mixing property.
%
%
\begin{lemme}
\label{dominationsm}
Let $t>0$ and $q>1$ be fixed. There exists a constant $A(t,q)$ such
that for any $x \in\Zd\setminus\{0\}$, for any
$A \in\mathcal F_t$, for any $B\in\mathcal F$, $\lambda\in\Lambda$
and every ${\ell} \ge1$,
\[
\bigl|\Pbarre_\lambda\bigl(A \cap(\tilde{\theta}_x^{\ell})^{-1}(B)\bigr)-\Pbarre
_\lambda(A) \Pbarre_{\T{\ell x}{\lambda}}(B)\bigr|\le A(t,q)q^{-{\ell}}.
\]
\end{lemme}
\begin{pf}
Let ${\ell} \ge1$. With Lemma \ref{supermagic}, we get
\begin{eqnarray*}
& & \bigl|\Pbarre_{\lambda}\bigl(A\cap\tilde{\theta}_x^{-\ell}(B)\bigr)-\Pbarre
_{\lambda}(A) \Pbarre_{\lambda}(\tilde{\theta}_x^{-\ell}(B))\bigr|\\
&&\qquad \le \bigl|\Pbarre_{\lambda}\bigl(t\le S_{\ell}(x),A\cap\tilde{\theta
}_x^{-\ell}(B)\bigr)-\Pbarre_{\lambda}\bigl(t\le S_{\ell}(x),A\bigr) \Pbarre
_{\lambda
}(\tilde{\theta}_x^{-\ell}(B))\bigr|\\
&&\qquad\quad{} +2\Pbarre
_{\lambda
}\bigl(t>S_{\ell}(x)\bigr)\\
&&\qquad = 2\Pbarre_{\lambda}\bigl(t>S_{\ell}(x)\bigr).
\end{eqnarray*}
Let us now fix $\alpha>0$.

With the Markov inequality, $\Pbarre_{\lambda}(S_{\ell}(x) \le t)
\le
\exp(\alpha t) \Ebarre_\lambda(\exp(-\alpha S_{\ell}(x)))$. Using the
two last points of Corollary \ref{invariancePbarre}, one has
\begin{eqnarray*}
\Ebarre_\lambda(\exp(-\alpha S_{\ell}(x)))
& \le&\Ebarre_\lambda\Biggl( \exp\Biggl( -\alpha\sum_{j=0}^{{\ell}-1}
\sigma(x) \circ\tilde{\theta}_x^j \Biggr) \Biggr) \\
& \le& \prod_{j=0}^{{\ell}-1}\Ebarre_\lambda\bigl( \exp\bigl(- \alpha
\sigma(x) \circ\tilde{\theta}_x^j\bigr) \bigr)\\
&=& \prod_{j=0}^{{\ell}-1}\Ebarre_{\T{jx}{\lambda}} ( \exp(- \alpha
\sigma(x))).
\end{eqnarray*}
Now we just have to prove the existence of some $\alpha>0$ such that
for every $\lambda\in\Lambda$,
\[
\Ebarre_{{\lambda}} ( \exp(- \alpha\sigma(x))) \le q^{-1}.\vadjust{\goodbreak}
\]
Let $\rho$ be the constant given in (\ref{uniftau}):
\begin{eqnarray*}
\Ebarre_\lambda( \exp(- \alpha\sigma(x))) &\le&\frac1\rho\E
_\lambda
(\exp(- \alpha\sigma(x))) \\
&\le&\frac1\rho\E_{\lambda_{\max
}}(\exp(-
\alpha\sigma(x)))\\
&\le&\frac1{\rho}\frac{2d\lambda_{\max}}{\alpha
+2d\lambda
_{\max}},
\end{eqnarray*}
because $\sigma(x)\ge t(x)$, and $t(x)$ stochastically dominates an
exponential random variable with parameter $2d\lambda_{\max}$. This
gives the desired inequality if $\alpha$ is large enough.
\end{pf}

We can now move forward to the proof of the ergodicity properties of
the systems $(\Omega,\mathcal{F},\Pbarre,\tilde{\theta}_x)$.
\begin{pf*}{Proof of Theorem \ref{systemeergodique}}
We have already seen in Corollary \ref{invariancePbarre} that for any
$x \in\Zd$, the probability measure $\Pbarre$ is invariant under the
action of~$\tilde{\theta}_x$. To prove ergodicity, we use an embedding
in a larger space to consider simultaneously a random environment and a
random contact process.\looseness=1

We thus set $\tilde{\Omega}=\Lambda\times\Omega$, equipped with the
$\sigma$-algebra $\tilde{\mathcal{F}}=\mathcal{B}(\Lambda)\otimes
\mathcal{F}$,
and we define a probability measure $\Qbarre$ on $\tilde{\mathcal
{F}}$ by
\[
\forall(A,B)\in\mathcal{B}(\Lambda)\times\mathcal{F} \qquad\Qbarre
(A\times
B)=\int_{\Lambda} \one_A(\lambda)\Pbarre_{\lambda}(B) \,d\nu
(\lambda).
\]
We define the transformation $\tilde{\Theta}_x$ on $\tilde{\Omega}$
by setting
$\tilde{\Theta}_x(\lambda,\omega)=(x.\lambda,\tilde{\theta
}_x(\omega
))$. It is easy to see that $\Qbarre$ is invariant under
$\tilde{\Theta}_x$. Indeed, for
$(A,B)\in\mathcal{B}(\Lambda)\times\mathcal{F}$, using Lemma \ref
{magic}, one has
\begin{eqnarray*}
\Qbarre\bigl(\tilde{\Theta}_x(\lambda, \omega)\in A\times B\bigr) & = &
\Qbarre
\bigl(x.\lambda\in A, \tilde{\theta}_x(\omega)\in B\bigr)\\
& = & \int_{\Lambda}\one_A(x.\lambda)\Pbarre_{\lambda}\bigl(\tilde
{\theta
}_x(\omega)\in B\bigr) \,d\nu(\lambda)\\
& = & \int_{\Lambda}\one_A(x.\lambda)\Pbarre_{x.\lambda}(B)
\,d\nu(\lambda
)\\
& = & \int_{\Lambda}\one_A(\lambda)\Pbarre_{\lambda}(B) \,d\nu
(\lambda)\\
&=& \Qbarre(A\times B).
\end{eqnarray*}
Note that if $g(\lambda,\omega)=f(\lambda)$, then $\int g \,d\Qbarre
=\int f \,d\nu$.\vspace*{1pt}

Similarly, if $g(\lambda,\omega)=f(\omega)$, then $\int g \,d\Qbarre
=\int f \,d\Pbarre$.

Note that $\mathcal{A}=\bigcup_{t\ge0}\mathcal{F}_t$ is an
algebra that generates $\mathcal{F}$. To prove that $\tilde{\theta}_x$
is ergodic, it is then sufficient to show that for every $A\in\mathcal{A}$,
%
%
\begin{equation}
\label{lacon}
\frac1{n}\sum_{k=0}^{n-1}\one_A(\tilde{\theta}_x^k)
\qquad\mbox{converges in } L^2(\Pbarre) \mbox{ to }\Pbarre(A).
\end{equation}
The quantity above can be seen as a function of the two
variables $(\lambda,\omega)$.
Thus it is equivalent to prove that the sequence of functions
$(\lambda,\omega)\mapsto\frac1{n}\sum_{k=0}^{n-1}\one_A(\tilde
{\theta}_x^k\omega)$ converges to $\Pbarre(A)$ in $L^2(\Qbarre)$.
Let $A\in\mathcal{A}$ and $t>0$ be such that $A\in\mathcal{F}_t$. For
every $(\omega,\lambda)\in\tilde{\Omega}$, we split the sum into two
terms:
\begin{eqnarray*}
\frac1{n}\sum_{k=0}^{n-1}\one_A(\tilde{\theta}_x^k\omega)
&=&
\frac1{n}\sum_{k=0}^{n-1}\bigl(\one_A(\tilde{\theta}_x^k\omega)-\Pbarre
_{kx.\lambda}(A)\bigr)\\
&&{}+\frac1{n}\sum_{k=0}^{n-1}\Pbarre
_{kx.\lambda}(A).
\end{eqnarray*}
If we set $f(\lambda)=\Pbarre_{\lambda}(A)$, the second term can be written
\[
\frac1{n}\sum_{k=0}^{n-1}\Pbarre_{kx.\lambda}(A)=\frac1{n}\sum
_{k=0}^{n-1} f(kx.\lambda).
\]
Since $x\in\Erg(\nu)$, the Von Neumann ergodic theorem says that this
quantity converges in $L^2(\nu)$ to $\int fd\nu=\Pbarre(A)$.
Seen as a function of $(\lambda,\omega)$, it also converges in
$L^2(\Qbarre)$ to $\Pbarre(A)$.
Set, $k\ge0$,
\[
Y_k=\one_A(\tilde{\theta}_x^k\omega)-\Pbarre_{\lambda}(\tilde
{\theta
}_x^{-k}(A))=\one_A(\tilde{\theta}_x^k\omega)-\Pbarre_{kx.\lambda}(A)
\]
%
and $L_n=Y_0+Y_1+\cdots+Y_{n-1}$.
It only\vspace*{1pt} remains to prove that $L_n/n$ converges to $0$ in $L^2(\Qbarre)$.
As $Y_k=Y_0\circ\tilde{\Theta}_x^k$, the field $(Y_k)_{k\ge0}$ is
stationary. We thus have
\begin{eqnarray*}
\int L_n^2 \,d\Qbarre& = & \sum_{0\le i,j\le n-1}\int Y_i Y_j
\,d\Qbarre\\
& =& \sum_{i=0}^{n-1} \int Y_i^2 \,d\Qbarre+2\sum_{\ell=1}^{n-1}
(n-\ell
) \int Y_0 Y_{\ell} \,d\Qbarre\\
& \le& 2n \Biggl(\sum_{\ell=0}^{+\infty} \biggl|\int Y_0 Y_{\ell}
\,d\Qbarre\biggr| \Biggr)
\le2n \Biggl(\sum_{\ell=0}^{+\infty} \int_{\Lambda
}|\Ebarre_{\lambda} (Y_0 Y_{\ell})| \,d\nu(\lambda) \Biggr)\\
& \le& 2n \Biggl(\sum_{\ell=0}^{+\infty} \int_{\Lambda}\bigl| \Pbarre
_{\lambda}\bigl(A\cap\tilde{\theta}_x^{-\ell}(A)\bigr)-\Pbarre_{\lambda}(A)
\Pbarre_{\lambda}(\tilde{\theta}_x^{-\ell}(A))\bigr| \,d\nu(\lambda)
\Biggr)\\
& \le& 2n \Biggl(\sum_{\ell=0}^{+\infty} A(t,2)2^{-\ell} \Biggr)=4A(t,2)n,
\end{eqnarray*}
thanks to Lemma \ref{dominationsm}. This ends the proof of (\ref
{lacon}), hence, the proof of Theorem \ref{systemeergodique}.
\end{pf*}

\section{Bound for the lack of subadditivity}
\label{controlediff}
In this section, we are going to bound quantities such as $\sigma
(x+y)-[\sigma(x)+\sigma(y)\circ\tilde{\theta}_x]$ and $\sigma(x)-t(x)$.

We will use these results in the application of a (almost) subadditive
ergodic theorem in Section \ref{forme}. In both cases, we use a kind of
restart argument. Considering the definition of the essential hitting
time $\sigma$, we will have to deal with two types of sums of random
variables that are quite different: sums of $v_i-u_i$ on one hand, and
sums of $u_{i+1}-v_i$ on the other hand.
\begin{itemize}
\item The life time $v_i(x) -u_i(x)$ of the contact process starting
from $x$ at time $u_i(x)$ can be bounded independently of the precise
configuration of the process at time $u_i(x)$. So the control is quite simple.
\item On the contrary, $u_{i+1}(x)-v_i(x)$, which represents the amount
of time needed to reinfect site $x$ after time $v_i(x)$, clearly
depends on the whole configuration of the process at time $v_i(x)$,
which is not easy to control precisely and uniformly in $x$. This
explains why the restart argument we use is more complex and the
estimates we obtain less accurate than in more classical situations
(e.g., in Section \ref{restartestimees}, we obtain the
exponential estimates of Proposition \ref{propuniforme} by standard
restart arguments).
\end{itemize}
As an illustration of the first point, we easily obtain the following lemma.
%
%
\begin{lemme}
\label{LEMtpsdevie}
There exist $A,B>0$ such that for every $\lambda\in\Lambda$,
%
%
\begin{equation}
\forall x \in\Zd,\forall t>0 \qquad\Pbarre_{\lambda}\bigl(\exists
i<K(x)\dvtx v_i(x)-u_i(x) >t\bigr) \le A\exp(-Bt).\hspace*{-20pt}
\end{equation}
\end{lemme}
\begin{pf}
Let $F\dvtx\Omega\to\R_+$ be a measurable function and $x \in\Zd$. We set
\[
\mathcal{L}_x(F)=\sum_{i=0}^{+\infty}\one_{\{u_i(x)<+\infty\}}
F\circ
\theta_{u_i(x)}.
\]
With the Markov property and the definition of $K(x)$, we have
\begin{eqnarray*}
\E_{\lambda}[\mathcal{L}_x(F)] & = &\sum_{i=0}^{+\infty}\E
_{\lambda}\bigl[\one
_{\{u_i(x)<+\infty\}}\bigr]\E_{\lambda}[ F]= \Biggl(1+\sum
_{i=0}^{+\infty}\P
_{\lambda}\bigl(K(x)>i\bigr)\Biggr)\E_{\lambda}[F]\\
& =& \bigl(1+\E_{\lambda}[K(x)]\bigr)\E_{\lambda}[F]\le\biggl( 1+\frac1\rho
\biggr)\E_{\lambda}[F],
\end{eqnarray*}
where the last equality comes from Lemma \ref{Kgeom}. We choose
$F=\one
_{\{t<u_i(x)-v_i(x)<+\infty\}}$, and with estimate (\ref{uniftau}),
we obtain
\begin{eqnarray*}
\Pbarre_{\lambda}\bigl(\exists i<K(x)\dvtx v_i(x)-u_i(x) >t\bigr) & \le&
\frac{1}{\rho}\P_{\lambda}\bigl(\exists i<K(x)\dvtx v_i(x)-u_i(x) >t\bigr) \\
& \le& \frac{1}{\rho} \P_{\lambda}\bigl(\mathcal{L}_x(F) \ge1\bigr)
\le\frac{1}{\rho} \E_{\lambda}[\mathcal{L}_x(F)] \\
& \le& \frac{1}{\rho}\biggl( 1+\frac1\rho\biggr) \P_{\lambda
}(t<\tau
^x<+\infty).
\end{eqnarray*}
We can then conclude with inequality (\ref{grosamasfinis}).
\end{pf}

To deal with the reinfection times $u_{i+1}(x)-v_i(x)$, the idea is to
look for a point $(y,t)$ (in space--time coordinates) close to
$(x,u_i(x))$, infected from $(0,0)$ and with infinite life time. The
at-least-linear-growth estimate (\ref{retouche}) will then ensure it
does not take too long to reinfect $x$ after time $t$, just by looking
at infection starting from the new source point $(y,t)$. The difficulty
lies in the control of the distance between $(x,u_i(x))$ and a source
point $(y,t)$; if the configuration around $(x,u_i(x))$ is
``reasonable,'' this point will not be too far from $(x,u_i(x))$, and we
will obtain a good control of $u_{i+1}(x)$ and~$u_i(x)$.

We recall that for every $x\in\Zd$, $\omega_x$ is the Poisson point
process giving the possible death times at site $x$, and that $M$ and
$c$ are, respectively, given in~(\ref{richard}) and (\ref{retouche}).
Note that we can assume that $M>1$. We note
%
%
\begin{equation}\label{gamma}
\gamma=3M(1+1/c)>3.
\end{equation}
For $x,y\in\Zd$ and $t>0$, we say that the growth from $(y,0)$ is bad
at scale $t$ with respect to $x$ if the following event occurs:
\begin{eqnarray*}
E^y(x,t)& =& \{\omega_y[0, t/2]=0\} \cup\{H^y_{ t}\not\subset y+B_{M
t}\} \cup \{ t/2<\tau^y<+\infty\} \\
&&{} \cup\bigl\{\tau^y=+\infty, \inf\{
s\ge2
t\dvtx x\in\xi^y_{s}\}> \gamma t\bigr\}.
\end{eqnarray*}
We want to check that with a high probability, there is no such bad
growth point in a box around $x$. So we define, for every $x\in\Zd$,
every $L> 0$ and every $t>0$,
\[
N_L(x,t)=\sum_{y\in x+B_{M t+2}} 
\int_0^L \one_{E^y(x,t)}\circ\theta_s \,d \biggl( \omega_y+\sum
_{e\in\Ed
\dvtx y\in e}\omega_e+\delta_0 \biggr)(s).
\]
In other words, we count the number of points $(y,s)$ in the space--time
box $(x+B_{M t+1})\times[0,L]$ such that something happens for site $y$
at time $t$, either a possible death, or a possible infection, and at
this time the bad event $E^y(x,t)\circ\theta_s$ occurs.
We first check that if the space--time box has no bad points and if
$u_i(x)$ is in the time window, then we can control the delay before
the next infection.
%
%
\begin{lemme}
\label{futfut}
If $N_L(x,t)\circ\theta_{s}=0$
and $s+ t\le u_i(x)\le s+L$, then
$v_i(x)=+\infty$ or $u_{i+1}(x)-u_i(x)\le\gamma t$.
\end{lemme}
\begin{pf}
By definition of $u_i(x)$, site $x$ is infected from $(0,0)$ at
time~$u_i(x)$. Since $s+t \le u_i(x)\le s+L$ and\vadjust{\goodbreak} $u_i(x)$ is a possible
infection time for $x$, the nonoccurrence of $E^x(x,t) \circ\theta
_{u_i(x)}$ ensures that $\tau^x\circ\theta_{u_i(x)}=+\infty$ or that
$\tau^x\circ\theta_{u_i(x)}\le t/2$.
If $\tau^x\circ\theta_{u_{i}(x)}=+\infty$, we are done because then
$v_{i}(x)=+\infty$.
Otherwise, note that $v_{i}(x)-u_i(x)\le t/2$.

By definition, there exists an infection path $\gamma_i\dvtx
[0,u_i(x)]\rightarrow\Z^d$ from $(0,0)$ to $(x,u_i(x))$, that is, such
that $\gamma_i(0)=0$ and $\gamma_i(u_i(x))=x$. Consider the portion of
$\gamma_i$ between time $u_i(x)- t$ and time $u_i(x)$. Denote by
$x_0=\gamma_i(u_i(x)-t)$ and let us see that $x_0 \in x+B_{M t+2}$.
Indeed, if $x_0 \notin x+B_{M t+2}$, we seek the first time $t_1$ after
time $u_i(x)- t$ when $\gamma_i$ enters in $x+B_{M t+2}$ at a site we
call $x_1$ (note that since $x_1$ is in the inside boundary of $x+B_{M
t+2}$, we have $\|x-x_1\|_\infty\ge Mt+1$). Time $t_1$ is a possible
infection time for $x_1$, and the nonoccurrence of $E^{x_1}(x,t) \circ
\theta_{t_1}$ ensures that the infection of $x$ from $(x_1,t_1)$ will
at least require a delay $t$, which contradicts $u_i(x)-t\ge0$.

So $x_0 \in x+B_{M t+2}$. Since $N_L(x,t)\circ\theta_{s}=0$, the first
possible death at site~$x_0$ after time $u_i(x)-t$ cannot occur after
a delay of $t/2$; thus the first time $t_2$ when the path $\gamma_i$
jumps to a different point $x_2$ satisfies $t_2 \le
u_i(x)-t+t/2=u_i(x)-t/2$. Consequently, when $(x_2,t_2)$ infects
$(x,u_i(x))$, it is at least $t/2$ aged, and the nonoccurrence of
$E^{x_2}(x,t) \circ\theta_{t_2}$ ensures it lives forever and
\[
\inf\{u\ge2 t\dvtx x \in\xi^{x_2}_{u}\}\circ\theta_{t_2}\le\gamma t.
\]
So there exists $t_3\in[t_2+2 t,t_2+ \gamma t]$ with $x\in\xi
^0_{t_3}$. Since $v_{i}(x)-u_i(x)\le t/2$,
one has
\[
t_3\ge t_2+2 t\ge\bigl(u_i(x)- t\bigr)+2 t=u_i(x)+ t\ge v_i(x).
\]
Finally, $u_{i+1}(x)-u_i(x)\le t_3-u_i(x)\le t_2-u_i(x)+ \gamma t\le
\gamma t$.
\end{pf}

Now we estimate the probability that a space--time box contains no bad points.
%
%
\begin{lemme}
\label{NL}
There exist $A_{\rref{eNL}},B_{\rref{eNL}}>0$ such that for every
$\lambda
\in\Lambda$,
%
%
\begin{equation}
\label{eNL}
\forall L>0, \forall x\in\Zd,\forall t>0\qquad \P_\lambda
\bigl(N_L(x,t)\ge1\bigr)\le
A_{\theequation}(1+L)\exp(-B_{\theequation}t).\hspace*{-35pt}
\end{equation}
\end{lemme}
\begin{pf}
Let us first prove there exist $A,B>0$ such that for every $\lambda\in
\Lambda$,
%
%
\begin{equation}
\label{eEE}\quad
\forall x\in\Zd,\forall t>0, \forall y\in x+B_{M t+2}\qquad
\P
_{\lambda}( E^y(x,t))\le A\exp(-Bt).
\end{equation}
Let $x \in\Zd$, $t>0$ and $y \in x+B_{M t+2}$.
If $\tau^y=+\infty$, there exists $z\in\xi^y_{2 t}$ with $\tau
^z\circ
\theta_{2t}=+\infty$.
Thus, definition (\ref{gamma}) of $\gamma$ implies that
\begin{eqnarray*}
& & \bigl\{\tau^y=+\infty, \inf\{s\ge2t\dvtx x\in\xi^y_{s}\}> \gamma
t\bigr\}
\\
&&\qquad \subset\{\xi^y_{2 t}\not\subset y+B_{2M t}\}
\cup\bigcup_{z\in y+B_{2M t}}
\{ t^z(x)\circ\theta_{2t}>(\gamma-2M) t \}\\
&&\qquad \subset \{\xi^y_{2t}\not\subset y+B_{2Mt}\}
\cup\bigcup_{z\in y+B_{2M t}}\biggl\{t^z(x) \circ\theta_{2t}>\frac
{\|
x-z\|}c+Mt-\frac{3}c \biggr\}.
\end{eqnarray*}
Hence, with (\ref{richard}) and (\ref{retouche}),
\begin{eqnarray*}
& &\P_{\lambda}(\tau^y=+\infty, \inf\{s\ge2 t\dvtx x\in\xi
^y_{s}\}>
\gamma t)\\
&&\qquad\le A\exp(-2BM t)+(1+4Mt)^d A\exp\bigl(-B(Mt-3/c)\bigr).
\end{eqnarray*}
The distribution of the number $\omega_y([0, t/2])$ of possible deaths
on site $y$
between time $0$ and time $t/2$ is a Poisson law with parameter $t/2$, so
\[
\P_{\lambda}\bigl(\omega_y([0, t/2])=0\bigr)= \exp(- t/2).
\]
The two remaining terms are controlled with (\ref{richard}) and (\ref
{grosamasfinis}); this gives (\ref{eEE}).

Now fix $y\in x+B_{Mt+2}$ and note $\beta_y=\omega_y+\sum
_{e\in\Ed\dvtx y\in e}\omega_e$.
Under $\P_{\lambda}$, $\beta_y$ is a Poisson point process with
intensity $2d\lambda_e$.
Let $S_0=0$ and $(S_n)_{n\ge1}$ be the increasing sequence of the
times given by this process:
\[
\int_0^L \one_{E^y(x,t)}\circ\theta_s \,d(\beta_y+\delta
_0)(s)=\sum
_{n=0}^{+\infty}\one_{\{S_n\le L\}}\one_{E^y(x,t)}\circ\theta_{S_n}.
\]
So, with the Markov property,
\begin{eqnarray*}
&& \E_{\lambda}\biggl(\int_0^L \one_{E^y(x,t)}\circ\theta_s
\,d(\beta
_y+\delta_0)(s)\biggr) \\
&&\qquad = \sum_{n=0}^{+\infty}\E_{\lambda}\bigl( \one_{\{S_n\le L\}
}\one
_{E^y(x,t)}\circ\theta_{S_n} \bigr)
= \sum_{n=0}^{+\infty}\E_{\lambda} \bigl( \one_{\{S_n\le L\}}
\bigr) \P
_{\lambda}(E^y(x,t))\\
&&\qquad = \bigl(1+ \E_{\lambda}[\beta_y([0,L])] \bigr) \P_{\lambda
}(E^y(x,t)) = \bigl(1+L(2d\lambda_e+1)\bigr)\P_{\lambda}(E^y(x,t)).
\end{eqnarray*}
So (\ref{eNL}) follows from (\ref{eEE}), from the remark that $\P
_\lambda(N_L(x,t) \ge1) \le\E_\lambda[N_L(x,t)]$ and from an obvious
bound on the cardinality of $B_{Mt+2}$.
\end{pf}

Once the process is initiated, Lemma \ref{futfut} can be used
recursively to control $u_{i+1}(x)-u_i(x)$. To initiate the process, we
assume that there exists a point $(u,s)$, reached from $(0,0)$, living
infinitely and close to $x$ in space.
%
%
\begin{lemme}
\label{futfut2}
For any $t,s>0$, for every $x \in\Zd$, the following inclusion holds:
%
%
\begin{eqnarray}
&&\{\tau^0=+\infty\}
\cap\{\exists u\in x+B_{Mt+2}, \tau_u\circ\theta_s=+\infty
,
u\in\xi^0_s\} \nonumber\\
\label{Ntamere}
&&\quad{}\cap\bigl\{N_{K(x)\gamma t}(x,t)\circ\theta_{s}=0\bigr\} 
\\ 
\label{pasdegrossaut}
&&\quad{} \cap\bigcap_{1\le i<K(x)} \{v_i(x)-u_i(x)< t\}
\\ 
\label{incluse}
&&\qquad\subset\{\tau^0=+\infty\} \cap\{\sigma
(x)\le s+K(x)\gamma t \}.
\end{eqnarray}
\end{lemme}
\begin{pf}
If every finite $u_i(x)$ is smaller than $s+t$, we are done because
$\sigma(x)\le s+t\le s+K(x)\gamma t$.
So set
\[
i_0=\max\{i\dvtx u_i(x)\le s+t\}.\vadjust{\goodbreak}
\]
Since $v_{i_0}(x)<+\infty$, the event (\ref{pasdegrossaut}) ensures that
$ v_{i_0}(x)-u_{i_0}(x)< t$, and so
$v_{i_0}(x)\le s+ 2t$.
Now, since $\tau^u=+\infty$, the nonoccurrence of $E^u(x,t) \circ
\theta
_s$ implied by (\ref{Ntamere})
says that
\[
\inf\{s \ge2t\dvtx x \in\xi^{u}_s\}\circ\theta_{s} \le\gamma t,
\]
which leads to $u_{i_0+1}(x) \le s+\gamma t$.
Noting that for any $j \ge1$,
$u_{i_0+j}(x+y)\ge s+ {t}$,
we prove by a recursive use of Lemma \ref{futfut} with the event $\{
N_{K(x)\gamma t}(x, t)\circ\theta_{s}=0\}$ that
\[
\forall j\in\{1,\ldots, K(x)-i_0\}\qquad u_{i_0+j}\le s+j\gamma{t}.
\]
For $j=K(x)-i_0$, we get
$
\sigma(x) = u_{i_0+j}(x)\le s+(K(x)-i_0)\gamma t
\le s+K(x)\gamma t,
$
which proves (\ref{incluse}).
\end{pf}

\subsection{Bound for the lack of subadditivity}

To bound $\sigma(x+y)-[\sigma(x)+\sigma(y)\circ\tilde{\theta
}_x]$, we
apply the strategy we have just explained around site $x+y$. To
initiate the recursive process, one can benefit here from the existence
of an infinite start at the precise point $(x+y, \sigma(x)+\sigma(y)
\circ\tilde{\theta}_y)$.
\begin{pf*}{Proof of Theorem \ref{presquesousadditif}}
Let $x,y \in\Zd$, $\lambda\in\Lambda$ and $t>0$. We set $s=\sigma
(x)+\sigma(y)\circ\tilde{\theta}_x$:
\begin{eqnarray*}
&& \Pbarre_{\lambda}\bigl( \sigma(x+y)>\sigma(x)+\sigma(y)\circ\tilde
{\theta
}_x+t\bigr) \\
&&\qquad \le \Pbarre_{\lambda}\biggl( K(x+y)> \frac{\sqrt t}{\gamma}
\biggr)\\
&&\qquad\quad{}+\Pbarre_{\lambda}
\biggl(\tau^0=+\infty, K(x+y) \le\frac{\sqrt t}{\gamma},
\sigma(x+y)\ge s + K(x+y)\gamma\sqrt t\biggr).
\end{eqnarray*}
With the sub-geometrical behavior of the tail of $K$ given in
Lemma \ref
{Kgeom} and the uniform control (\ref{uniftau}), we can control the
first term. Note that if $K(x+y) \le\frac{\sqrt t}{\gamma}$, then
$K(x+y) \gamma\sqrt t \le t$, and so that
\[
\bigl\{N_{K(x+y) \gamma\sqrt t}\bigl(x+y,\sqrt t\bigr) \ge1\bigr\} \subset\bigl\{
N_{t}\bigl(x+y,\sqrt t\bigr) \ge1\bigr\}.
\]
We apply Lemma \ref{futfut2} around $x+y$, on a scale $\sqrt t$, an
initial time $s=\sigma(x)+\sigma(y)\circ\tilde{\theta}_x$ and a source
point $u=x+y$,
%
%
\begin{eqnarray}\label{Clui}
&& \Pbarre_{\lambda} \biggl(
\tau^0=+\infty, K(x+y) \le\frac{\sqrt t}{\gamma}, \sigma
(x+y)\ge
s + K(x+y)\gamma\sqrt t \biggr) \nonumber\\
&&\qquad \le \Pbarre_{\lambda}\bigl(N_t\bigl(x+y, \sqrt t\bigr)\circ\theta_{s} \ge1\bigr)
\\
&&\qquad\quad{} +\Pbarre_\lambda\bigl(\exists i<K(x+y)\dvtx v_i(x+y)-u_i(x+y) >\sqrt t\bigr).
\nonumber
\end{eqnarray}
Since $N_t(x+y,\sqrt{t})=N_t(0,\sqrt{t})\circ T_x\circ T_y$ and
$s=\sigma(x)+\sigma(y)\circ\tilde{\theta}_x$, we have
\[
N_t\bigl(x+y,\sqrt t\bigr)\circ\theta_{s}=N_t\bigl(0,\sqrt{t}\bigr)\circ\tilde{\theta
}_y\circ\tilde{\theta}_x.
\]
Thus
$ \Pbarre_{\lambda}(N_t(x+y,\sqrt t)\circ\theta_s \ge1)= \Pbarre
_{(x+y).\lambda} (N_t(0,\sqrt t)\ge1)$, which is controlled by
Lemma \ref{NL} and estimate (\ref{uniftau}).
Finally, (\ref{Clui}) is bounded with Lemma \ref{LEMtpsdevie}.
\end{pf*}
%
%
\begin{coro}
\label{momentsecart}
For $x,y \in\Zd$, set $r(x,y)=(\sigma(x+y)-(\sigma(x)+\sigma
(y)\circ
\tilde{\theta}_x))^+$.

For any $p\ge1$, there exists $M_p>0$ such that
%
%
\begin{equation} \label{momecart}
\forall\lambda\in\Lambda,\forall x,y\in\Zd\qquad\Ebarre
_{\lambda
}[r(x,y)^p]\le M_p.
\end{equation}
\end{coro}
\begin{pf}
We write $\Ebarre_{\lambda}[r(x,y)^p]=\int_0^{+\infty}
pu^{p-1}\Pbarre
_{\lambda}(r(x,y)>u) \,du$ and use Theorem \ref{presquesousadditif}.
\end{pf}

\subsection{Control of the discrepancy between hitting times and
essential hitting times}

To bound $\sigma(x)-t(x)$, we would like to apply the same strategy
starting from $(x,t(x))$ but we do not have any natural candidate for
an infinite start close to this point. We are going to look for such a
point along the infection path between $(0,0)$ and $(x,t(x))$ which
requires controls on a~space--time box whose height (in time) of order
$t(x)$, that is, of order $\|x\|$. So we will lose in the precision of
the estimates and in their uniformity.
%
%
\begin{prop}
\label{importante}
There exist $A_{\rref{eimportante}},B_{\rref{eimportante}},\alpha
_{\rref
{eimportante}}>0$ such that for every \mbox{$z>0$}, every $x\in\Zd$, every $
\lambda\in\Lambda$,
%
%
\begin{equation} \label{eimportante}\quad
\Pbarre_\lambda\bigl(\sigma(x)\ge t(x)+K(x)\bigl(\alpha_{\theequation
}\log
(1+\|x\|) +z\bigr)\bigr)
\le A_{\theequation}\exp(-B_{\theequation}z).
\end{equation}
\end{prop}
\begin{pf} For $x,y\in\Zd$ and $t,L>0$, we define
\begin{eqnarray*}
\tilde{E}^y(t) & = &\biggl\{\tau_y<+\infty, \bigcup_{s\ge0}H^y_s\not
\subset y+B_{Mt}\biggr\}, \\
\tilde N_L(x,t) & = & \sum_{y\in x+B_{Mt+1}}\int_0^L \one_{\tilde
E^y(t)}\circ\theta_s d \biggl(\sum_{e\in\Ed\dvtx y\in e}\omega_e
\biggr)(s).
\end{eqnarray*}
With (\ref{uniftau}), (\ref{richard}) and (\ref{grosamasfinis}), it is
easy to get the existence of $A,B>0$ such that
%
%
\begin{equation}\label{hauteur}
\forall\lambda\in\Lambda,\forall x \in\Zd,\forall t>0\qquad
\Pbarre_{\lambda}\bigl(\tilde N_L(x,t)\ge1\bigr)\le A (1+L)\exp(-B
t).\hspace*{-25pt}
\end{equation}
Now, we choose the last point $(u,s)$ on the infection path between
$(0,0)$ and $(x,t(x))$ such that
$\tau^u \circ\theta_s=+\infty$. Note that on $\{\tau^0=+\infty\}$,
such an $s$ always exists.

Let us see that if $\tilde N_{t(x)}(x,t)=0$, then $u\in x+B_{Mt+2}$.
Indeed, if $\|u-x\|> Mt +2$, we consider the first point $(u',s')$ on
the infection path after $(u,s)$ to be in $x+B_{Mt}$. The definition of
$s$ ensures that the contact process starting from $(u',s')$ does not
survive, but, since it contains $(x,t(x))$, its diameter must be larger
than $Mt$, which implies that $\tilde N_{t(x)}(x,t)\ge1$, and gives
the desired implication.

On event $\{\tilde N_{t(x)}(x,t)=0\}$, we are going to apply Lemma \ref
{futfut2} around point $(x,0)$, at scale
\[
t=\frac{\alpha\log(1+\|x\|)+z}{\gamma}\ge\frac{z}{\gamma}
\]
with source point
$(u,s)$ and a time length $L=K(x)\gamma t$. Here and in the following,
$\alpha>0$ is a large constant that will be chosen later. Since $s\le t(x)$,
%
%
\begin{eqnarray} \label{MAJhic}
&& \Pbarre_\lambda\bigl(\sigma(x)\ge t(x)+K(x)\bigl(\alpha\log(1+\|x\|)
+z\bigr)\bigr)\nonumber\\
&&\qquad= \Pbarre_\lambda\bigl(\sigma(x)\ge t(x)+K(x)\gamma t \bigr)
\nonumber\\
&&\qquad \le \Pbarre_\lambda\bigl(\sigma(x)\ge s+K(x)\gamma t \bigr)
\nonumber\\
&&\qquad \le \Pbarre_\lambda\bigl(\sigma(x)\ge s+K(x)\gamma t, \tilde
N_{t(x)}(x,t)=0 \bigr)\nonumber\\[-8pt]\\[-8pt]
&&\qquad\quad{} +\Pbarre_\lambda\bigl(\tilde
N_{t(x)}(x,t)\ge1
\bigr) \nonumber\\
&&\qquad \le \Pbarre_{\lambda}\bigl(N_{K(x)\gamma t}(x, t)\circ\theta_{s} \ge1\bigr)
\nonumber\\
&&\qquad\quad{}+\Pbarre_\lambda\bigl(\exists i<K(x)\dvtx v_i(x)-u_i(x)
>t\bigr)\nonumber\\
&&\qquad\quad{} + \Pbarre_\lambda\bigl(\tilde N_{t(x)}(x,t)\ge1 \bigr). \nonumber
\end{eqnarray}
The second term in (\ref{MAJhic}) is bounded with Lemma \ref
{LEMtpsdevie}. For the last term, we write
\[
\Pbarre_\lambda\bigl(\tilde N_{t(x)}(x,t)\ge1 \bigr) \le
\Pbarre
_\lambda\bigl(\tilde N_{{\|x\|}/c+z}(x,t)\ge1 \bigr) +\Pbarre
_\lambda\biggl( t(x)> \frac{\|x\|}c+z \biggr).
\]
The second term is controlled with (\ref{uniftau}) and (\ref
{retouche}), and (\ref{hauteur}) ensures that
\begin{eqnarray*}
\Pbarre_\lambda\bigl(\tilde N_{{\|x\|}/c+t}(x,t)\ge1 \bigr)
& \le& A \biggl( 1+\frac{\|x\|}c+z \biggr)\exp(-B t) \\
& \le& A \biggl( 1+\frac{\|x\|}c+z \biggr)\exp\biggl(-\frac
{B(\alpha
\log(1+\|x\|)+z)}{\gamma} \biggr) \\
& \le& A' \exp(-B'z)
\end{eqnarray*}
as soon as $\alpha$ is large enough.

For the first term of (\ref{MAJhic}), we note that
$N_{K(x)\gamma t}(x,t)\circ\theta_{s}\le N_{t(x)+K(x)\gamma t}(x,t)$.
Thus
\begin{eqnarray*}
\Pbarre_{\lambda}\bigl(N_{K(x)\gamma t}(x,t)\circ\theta_{s}\ge1\bigr)
&\le&\Pbarre_{\lambda} \bigl(N_{{\|x\|}/c+z+K(x)\gamma
t}(x,t)\ge
1\bigr)\\
&&{}+\Pbarre_{\lambda}\biggl(t(x)\ge\frac{\|x\|}c+z\biggr).
\end{eqnarray*}
As previously stated, the second term is bounded with (\ref{uniftau})
and (\ref{retouche}) and while using~(\ref{eNL}), we get
\begin{eqnarray*}
&& \Pbarre_{\lambda} \bigl(N_{{\|x\|}/c+z+K(x)\gamma t}(x,t)\ge
1\bigr) \\[-1pt]
&&\qquad\le \sum_{k=1}^{+\infty}\sqrt{\Pbarre_{\lambda}\bigl(K(x)=k\bigr)}\sqrt
{\Pbarre_{\lambda} \bigl(N_{k\gamma t+{\|x\|}/c+z}(x,t)\ge
1
\bigr)}\\[-1pt]
&&\qquad\le \sum_{k=1}^{+\infty}\sqrt{\Pbarre_{\lambda}\bigl(K(x)=k\bigr)}\sqrt
{A_{\rref{eNL}} \biggl(k\gamma t+\frac{\|x\|}c+z \biggr)\exp
(-B_{\rref
{eNL}}t)}\\[-1pt]
&&\qquad \le \sqrt{ A_{\rref{eNL}} \biggl(1+\frac{\|x\|}c \biggr)
(1+z)(1+\gamma t) }\exp\biggl(-\frac{B_{\rref{eNL}}}2 t \biggr)\\[-1pt]
&&\qquad\quad{}\times\sum
_{k=1}^{+\infty}\sqrt{(1+k)\Pbarre_{\lambda}\bigl(K(x)=k\bigr)}.
\end{eqnarray*}
The sub-geometrical behavior of the tail of $K(x)$ given by Lemma
\ref{Kgeom} ensures that the sum is finite, and we end the proof by
increasing $\alpha$ if necessary.\vspace*{-2pt}
\end{pf}
%
%
\begin{lemme}
For every $p \ge1$, there exists $C_{\rref{momv}}(p)>0$ such that for
every $x \in\Zd$
%
%
\begin{equation} \label{momv}
\forall\lambda\in\Lambda\qquad\Ebarre_{\lambda} \bigl(|\sigma(x)-t(x)|^p\bigr)
\le C_{\theequation}(p) \bigl(\log(1+\|x\|)\bigr)^{p}.\vspace*{-2pt}
\end{equation}
\end{lemme}
\begin{pf}
Set\vspace*{1pt} $V_x=\frac{\sigma(x)-t(x)}{K(x)}-\alpha_{\rref{eimportante}}\log
(1+\|
x\|)$.
By Proposition \ref{importante}, there exists a random variable $W$
with exponential moments that stochastically dominates $V_x$ under
$\Pbarre_{\lambda}$ for every $x$ and every $\lambda$. Moreover,
Lemma \ref{Kgeom} ensures that $K(x)$ is stochastically dominated by a
geometrical random variable $K'$.

Set $v(x)=\sigma(x)-t(x)=K(x)(\alpha\log(1+\|x\|)+V_x)$ and let
$p\ge
1$. With the Minkowski inequality, we have
\begin{eqnarray*}
(\Ebarre_{\lambda} v(x)^p)^{1/p}& \le& \alpha\log(1+\|x\|)(\Ebarre
_{\lambda} K(x)^p)^{1/p}+(\Ebarre_{\lambda} [ K(x)^p V_x^p
]
)^{1/p}\\
& \le& \alpha\log(1+\|x\|)(\Ebarre_{\lambda} K(x)^p)^{1/p}+
(\Ebarre_{\lambda} K(x)^{2p} \Ebarre_{\lambda}V_x^{2p}
)^{1/({2p})}\\
& \le& \alpha\log(1+\|x\|)(\E K'^p)^{1/p}+(\E K'^{2p} \E
W^{2p})^{1/({2p})},
\end{eqnarray*}
and the proof is complete.\vspace*{-2pt}
\end{pf}
%
%
\begin{coro}
\label{pareil}
$\Pbarre$-a.s., $ 
\lim_{\|x\|\to+\infty} \frac{|\sigma(x)-t(x)|}{\|x\|}=0$.\vspace*{-2pt}
\end{coro}
\begin{pf}
Let $p> d$. Equation (\ref{momv}) gives
\[
\sum_{x\in\Zd} \Ebarre\frac{|\sigma(x)-t(x)|^p}{(1+\|x\|)^p}
\le
C_{\rref{momv}}(p) \sum_{x\in\Zd} \frac{(\log(1+\|x|))^{p}}{(1+\|
x\|
)^p}<+\infty.
\]
So $(\frac{|\sigma(x)-t(x)|}{(1+\|x\|)})_{x \in\Zd}$ is
almost surely in $\ell^p(\Zd)$ and thus goes to $0$.\vadjust{\goodbreak}
\end{pf}
%
%
\begin{coro}
\label{sigmaa}
There exist $A_{\rref{esigmaa}},B_{\rref{esigmaa}},C_{\rref{esigmaa}}>0$
such that for every $\lambda\in\Lambda$,
%
%
\begin{equation}\label{esigmaa}\quad
\forall x \in\Zd,\forall t>0\qquad \Pbarre_\lambda
\bigl(\sigma
(x)\ge C_{\theequation}\|x\|+t\bigr) \le A_{\theequation}\exp
\bigl(-B_{\theequation}\sqrt{t}\bigr).
\end{equation}
\end{coro}
\begin{pf}
Let $\alpha=\alpha_{\rref{eimportante}}$ be given as in
Proposition \ref
{importante}, and note that if
$K(x)\le\frac1{2\alpha} \sqrt{\|x\|+t/2}$ and $z=\alpha\sqrt{\|x\|+t/2}$,
then, since $\log(1+u)\le\sqrt u$, we get
\[
K(x)[\alpha\log(1+\|x\|)+z]\le2zK(x)\le\|x\|+t/2.
\]
Thus with (\ref{retouche}) and (\ref{uniftau}),
\begin{eqnarray*}
&& \Pbarre_\lambda\biggl(\sigma(x)> \biggl(\frac1{c}+1\biggr)\|x\|
+t
\biggr) \\
&&\qquad \le \Pbarre_\lambda\biggl(t(x)\ge\frac{\|x\|}{c}+t/2 \biggr) +
\Pbarre_\lambda\biggl( K(x)>\frac1{2\alpha} \sqrt{\|x\|+t/2}
\biggr)\\
&&\qquad\quad{} + \Pbarre_\lambda\bigl(\sigma(x)>
t(x)+K(x)\bigl(\alpha
\log(1+\|x\|) +\alpha\sqrt{\|x\|+t/2}\bigr)\bigr).
\end{eqnarray*}
The first term is controlled with (\ref{retouche}), the second one with
Lemma \ref{Kgeom} and the last one by Proposition \ref{importante}.
\end{pf}
%
%
\begin{coro}
\label{propmoments}
For any $p\ge1$, there exists $C_{\rref{moms}}(p)>0$ such that
%
%
\begin{equation}
\label{moms}
\forall\lambda\in\Lambda,\forall x\in\Zd\qquad\Ebarre
_{\lambda
}[\sigma(x)^p] \le C_{\theequation}(p) (1+\|x\|)^{p}.
\end{equation}
\end{coro}
\begin{pf}
With the Minkowski inequality, one has
\[
(\Ebarre_{\lambda}[\sigma(x)^p])^{1/p}\le C_{\rref{esigmaa}}\|x\|
+\bigl(\Ebarre_{\lambda}\bigl[\bigl(\bigl(\sigma(x)-C_{\rref{esigmaa}}\|x\|\bigr)^+\bigr)^p\bigr]\bigr)^{1/p}.
\]
Moreover,
\[
\Ebarre_{\lambda}\bigl[\bigl(\bigl(\sigma(x)-C_{\rref{esigmaa}}\|x\|\bigr)^+\bigr)^p\bigr]=
\int_0^{+\infty} pu^{p-1}\Pbarre_{\lambda}\bigl(\sigma(x)-C_{\rref
{esigmaa}}\| x\|>u\bigr) \,du<+\infty
\]
by Corollary \ref{sigmaa}.
\end{pf}
\begin{Remark*}
In classical restart arguments, the existence of exponential moments
for a random variable usually comes from the following argument: if
$(X_n)_{n \in\N}$ are independent identically distributed random
variables with exponential moments, if $K$ is independent of the
$(X_n)_{n \in\N}$'s and also has exponential moments, then $
\sum_{0\le n \le K} X_n$ has exponential moments. Here, our
difficulties to precisely bound the reinfection times $u_{i+1}-v_i$
prevent us to use this scheme; we thus have to use {ad hoc} arguments,
which lead to weaker estimates.
\end{Remark*}

\section{Asymptotic shape theorems}
\label{forme}
We can now move forward to the proof of Theorem \ref{thFA}.
The first step consists of proving convergence for ratios of the type
$\frac{\sigma(nx)}n$.
With Corollary \ref{momentsecart}, we know that for every $n,p \ge0$,
\[
\Ebarre\bigl[\sigma\bigl((n+p)x\bigr)\bigr]\le\Ebarre[\sigma(nx)]+\Ebarre[\sigma(px)]+M_1.\vadjust{\goodbreak}
\]
Thus the Fekete lemma says that $\frac1{n}\Ebarre[\sigma(nx)]$ has a
finite limit when $n$ goes to $+\infty$
and the natural candidate for the limit of $\frac{\sigma(nx)}n$ is thus
\[
\mu(x)=\lim_{n\to+\infty} \frac{{\Ebarre} (\sigma(nx))}{n}.\vspace*{-2pt}
\]

%
\begin{theorem}
\label{cvdir}
$\Pbarre\mbox{-a.s. } \forall x \in\Zd\lim_{n \to
+\infty} \frac{\sigma(nx)}{n}= \lim_{n \to+\infty} \frac{\Ebarre
\sigma
(nx)}{n}=\mu(x)$.

This convergence also holds in any $L^p(\Pbarre)$, $p\ge1$.\vspace*{-2pt}
\end{theorem}

To prove this result, we need the two following (almost) subadditive
ergodic theorems, whose proof will be given in the \hyperref[app]{Appendix}.\vspace*{-2pt}
%
%
\begin{theorem}
\label{therg}
Let $(\Omega,\mathcal{F},\P)$ be a probability space, $(\theta
_n)_{n\ge
1}$ a collection of transformations leaving the probability measure $\P
$ invariant. On this space, we consider a collection $(f_n)_{n\ge1}$
of integrable functions, a collection
$(g_n)_{n\ge1}$ of nonnegative functions and a collection
$(r_{n,p})_{n,p\ge1}$ of real functions such that
%
%
\begin{equation}
\label{maiscestimportant}
\forall n,p\ge1\qquad f_{n+p}\le f_{n}+f_{p}\circ\theta_n+g_p\circ
\theta
_n+r_{n,p}.
\end{equation}
We assume that:
\begin{itemize}
\item$c=\inf_{n\ge1} \frac{\E f_n}n>-\infty$.
\item$g_1$ is integrable, $g_n/n$ almost surely converges to $0$ and
$\frac{\E g_n}n$ converges to~$0$.
\item There exists $\alpha>1$ and a sequence of positive numbers
$(C_p)_{p\ge1}$ such that
$\E[ (r^+_{n,p})^{\alpha}]\le C_p$ for every $n,p$ and
\[
\sum_{p=1}^{+\infty} \frac{C_p}{p^{\alpha}}<+\infty.
\]
\end{itemize}
Then $\frac1{n}\E f_n$ converges; if $\mu$ denotes its limit, one has
\[
\E\biggl[\liminf_{n \to+\infty} \frac{f_n}{n}\biggr]\ge\mu.
\]
If we set $\underline{f}=\liminf_{n \to+\infty} \frac
{f_n}{n}$, then $\underline{f}$ is invariant under the action of each
$\theta_n$.\vspace*{-2pt}
\end{theorem}
%
%
\begin{theorem}
\label{thergdeux}
We keep the setting and assumptions of Theorem \ref{therg}.
We assume, moreover, that for every $k$,
\[
\frac1{n}\Biggl(f_{nk}-\sum_{i=0}^{n-1}f_k\circ(\theta_k)^i
\Biggr)^+\to
0 \qquad\mbox{a.s.}
\]
Then $f_n/n$ converges a.s. to $\underline{f}$.\vspace*{-2pt}
\end{theorem}
\begin{pf*}{Proof of Theorem \ref{cvdir}}
We apply Theorem \ref{therg} with the choices $f_n=\sigma(nx)$,
$\theta
_n=\tilde{\theta}_{nx}$, $g_p=0$, $r_{n,p}=r(nx,px)$ and the
probability measure $\P=\Pbarre$. We take $\alpha>1$.
Corollary \ref{propmoments} gives the integrability of $\sigma(x)$
under $\Pbarre$ and Corollary~\ref{momentsecart} gives the necessary
controls on its moments.\vadjust{\goodbreak}

We now check the extra assumption of Theorem \ref{thergdeux}; it is
easy to see that
\[
t(nkx)\le\sum_{i=0}^{n-1}\sigma(kx)\circ(\tilde{\theta}_{kx})^i,\vspace*{-2pt}
\]
which implies that
$ (\sigma(nkx)-\sum_{i=0}^{n-1}\sigma(kx)\circ
(\tilde
{\theta}_{kx})^i)^+\le\sigma(nkx)-t(nkx)$.

Corollary \ref{pareil} ensures\vspace*{1pt} that this quantity is $o(n)$.
Thus $\sigma(nx)/n$ converges to a random variable $\mu(x)$, which is
invariant under the action of $\tilde{\theta}_x$. But Theorem \ref
{systemeergodique} says that this $\mu(x)$ is in fact a constant, which
ends the proof of the a.s. convergence.

To prove that a sequence converges in $L^p$, it suffices to show that
it converges a.s. and that it is bounded in $L^q$ for some $q>p$. Since
Corollary~\ref{propmoments} says that $f_n/n$ is bounded in any $L^p$,
the proof is complete.\vspace*{-2pt}
\end{pf*}

The next step is to prove the asymptotic shape result, namely,
Theorem~\ref{thFA}. We start by proving the shape result for the
essential hitting time $\sigma$, by following the classical strategy:
\begin{itemize}
\item We extend $\mu$ to an asymmetric norm on $\Rd$ in Lemma \ref{munorme}.
\item We prove that the directional convergence given by Theorem \ref
{cvdir} is in fact uniform in the direction in Lemma \ref{convunifsigma}.
\item We easily deduce the shape result from this lemma in Lemma \ref
{formesigma}.
\end{itemize}
To transpose this shape result for the classical hitting time $t$
(Lemma \ref{etpourt}), we just need to control the discrepancy between
$\sigma$ and $t$; this was done in Lemma \ref{pareil}. Finally, the
shape result for the coupled zone is proved in Lemma~\ref{zoneC} by
introducing a coupling time $t'$ and by bounding the difference between
this time $t'$ and the essential hitting time $\sigma$.

Note that we did not succeed in proving immediately that
$\mu$ could be extended to a norm, but only to an asymmetric norm; that
is, the property $\mu(\lambda x)=|\lambda|\mu(x)$
{a priori} only holds for nonnegative $\lambda$.
We will finally deduce from the asymptotic shape theorems that $\mu$ is
actually a norm.\vspace*{-2pt}
%
%
\begin{lemme}
\label{munorme}
$\!\!\!$The functional $\mu$ can be extended to an asymmetric norm~on~$\Rd$.\vspace*{-2pt}
\end{lemme}
\begin{pf}
\textit{Homogeneity in natural integers}.
By extracting subsequences, we prove the homogeneity in natural integers,
\[
\forall k\in\N,\forall x\in\Zd\qquad\mu(kx)=k\mu(x).
\]

\textit{Subadditivity}.
One has
$\sigma(nx+ny)\le\sigma(nx)+\sigma(ny)\circ\tilde{\theta
}_{nx}+r(nx,ny)$.

Since $\Pbarre$ is invariant under the action of $\tilde{\theta}_{nx}$,
we get, with Corollary \ref{momentsecart},
\[
\Ebarre\sigma(nx+ny) \le\Ebarre\sigma(nx)+\Ebarre\sigma
(ny)+\Ebarre
r(nx,ny)\le\Ebarre\sigma(nx)+\Ebarre\sigma(ny)+M_1.
\]
We deduce that
$\forall x \in\Zd,\forall y \in\Zd,\mu(x+y)\le\mu
(x)+\mu(y)$.

\textit{Extension to $\Rd$}.
The Fekete lemma ensures that
\[
\mu(x)+M_1=\inf_{n\ge1}\frac{\Ebarre{\sigma(nx)+M_1}}n,\vadjust{\goodbreak}
\]
so $\mu(x)\le\Ebarre\sigma(x)$.
Corollary \ref{propmoments} gives some $L>0$ such that $\Ebarre\sigma
(x)\le L\|x\|$ for any $x$. Finally, $\mu(x)\le L\|x\|$ for every $x
\in\Zd$, which leads to $|\mu(x)-\mu(y)|\le L\|x-y\|$. We can then extend
$\mu$ to $\Q^d$ par homogeneity, then to $\Rd$ by uniform continuity.

\textit{Positivity}.
Let $M$ be given by Proposition \ref{propuniforme}. With (\ref
{richard}), we obtain
\begin{eqnarray*}
\Pbarre\biggl(\sigma(nx)< \frac{n\|x\|}{2M}\biggr) & \le&
\Pbarre\biggl(t(nx)< \frac{n\|x\|}{2M}\biggr) \le\Pbarre\bigl(\xi
_{{n\|x\|}/{M}}\not\subset B_{{n\|x\|}/2} \bigr)\\
& \le&\int\frac{\P_{\lambda}(\tau^0=+\infty, \xi
^0_{{n\|x\|
}/({2M})}\not\subset B_{{n\|x\|}/2})}{\P_{\lambda}(\tau
^0=+\infty
)} \,d\nu(\lambda)\\
& \le&\frac{A}{\rho}\exp\biggl(-B\frac{n\|x\|}{2M}\biggr).
\end{eqnarray*}
With the Borel--Cantelli lemma, we deduce that $\mu(x)\ge\frac
{1}{2M}\|
x\|$.
This inequality, once established for every $x \in\Zd$, can be
extended by homogeneity and continuity to $\Rd$. So $\mu$ is an
asymmetric norm.
\end{pf}

In the following, we set $C=2C_{\rref{esigmaa}}$, where $C_{\rref
{esigmaa}}$ is as given in Corollary \ref{sigmaa}.
%
%
\begin{lemme}
\label{reglo}
For every $\epsilon>0$, $\Pbarre$-a.s., there exists $R>0$ such that
\[
\forall x,y \in\Zd\qquad(\|x\|\ge R\mbox{ and }\|x-y\|\le\epsilon\|
x\|)\Longrightarrow\bigl(|\sigma(x)-\sigma(y)|\le C\epsilon\|x\|\bigr).
\]
\end{lemme}
\begin{pf}
For $m \in\N$ and $\varepsilon>0$, we define the event
\[
A_m(\epsilon)=\{\exists x,y \in\Zd\dvtx\|x\|=m, \|x-y\|\le
\epsilon m \mbox{ and } |\sigma(x)-\sigma(y)|>C\epsilon m \}.
\]
Noting that
\[
A_m(\epsilon) \subset
\mathop{\bigcup_{(1-\epsilon) m\le\|x\| \le(1+\epsilon)m}}_{\|
x-y\|
\le\epsilon m}
\{\sigma(y-x)\circ\tilde{\theta}_{x}+r(x,y-x)>C\epsilon
m\},
\]
%
we see, with Corollaries \ref{sigmaa} and \ref{momentsecart}, that
%
\begin{eqnarray*}
\Pbarre_{\lambda} ( A_m(\epsilon) )
& \le&
\mathop{\sum_{(1-\epsilon) m\le\|x\| \le(1+\epsilon)m}}_{\|z\|\le
\epsilon m} \Pbarre_{{\lambda}}\bigl(\sigma(z)\circ\tilde{\theta
}_{x}+r(x,z)>C\epsilon m\bigr)\\
& \le&
\mathop{\sum_{(1-\epsilon) m\le\|x\| \le(1+\epsilon)m}}_{\|z\|\le
\epsilon m} \Pbarre_{\T{x}{\lambda}}\bigl(\sigma(z)>2C\epsilon m/3\bigr)\\
& &{} +\Pbarre_{{\lambda}}\bigl(r(x,y-x)>C\epsilon m/3\bigr)\\
& \le& (1+2\epsilon m)^d\bigl(1+2(1+\epsilon) m\bigr)^dA_{\rref{esigmaa}}\exp
\bigl(-B_{\rref{esigmaa}}\sqrt{C\epsilon m /3}\bigr)\\
& &{} +A_{\rref{momecart}}\exp\bigl(-B_{\rref{momecart}}\sqrt
{C'\epsilon m /3}\bigr)
\end{eqnarray*}
by Corollary \ref{sigmaa} and Theorem \ref{presquesousadditif}.
Integrating then with respect to $\lambda$, we conclude the proof with
the Borel--Cantelli lemma.
\end{pf}
%
%
\begin{lemme}
\label{convunifsigma}
$\Pbarre$-a.s. $ \lim_{\|x\|\to+\infty} \frac{|\sigma
(x)-\mu(x)|}{\|x\|}=0$.
\end{lemme}
\begin{pf}
Assume by contradiction that there exists $\epsilon>0$ such that the
event ``$|\sigma(x)-\mu(x)|> \epsilon\|x\|$ for infinitely many values
of $x$'' has a positive probability. We focus on this event. There
exists a random sequence $(y_n)_{n \ge0}$ of sites in $\Z^d$ such that
$\|y_n\|_1 \rightarrow+ \infty$ and, for every $n$, $|\sigma
(y_n)-\mu
(y_n)| \geq\varepsilon\|y_n\|_1$. By extracting a subsequence, we can
assume that
\[
\frac{y_n}{\|y_n\|_1} \rightarrow z.
\]
Fix $\varepsilon_1>0 $ (to be chosen later); we can find $z' \in\Erg
(\nu
)$ such that
\[
\biggl\| \frac{z'}{\|z'\|_1} -z\biggr\|_1\leq\varepsilon_1.
\]
For each $y_n$, we can find an integer point on $\R z'$ close to $y_n$.
Let $h_n$ be the integer part of $\frac{\|y_n\|_1}{\|z'\|_1}$. We have
\begin{eqnarray*}
\|y_n-h_n.z'\|_1 & \leq& \biggl\|y_n- \frac{\|y_n\|_1}{\|z'\|_1}
z'\biggr\|_1+ \biggl| \frac{\|y_n\|_1}{\|z'\|_1}- h_n\biggr| \|z'\|
_1 \\
& \leq&\|y_n\|_1 \biggl\|\frac{y_n}{\|y_n\|_1}-\frac{z'}{\|z'\|
_1}
\biggr\|_1 +\|z'\|_1.
\end{eqnarray*}
Take $N>0$ large enough to have
$(n \geq N) \Rightarrow(\|\frac{y_n}{\|y_n\|_1}-z\|
_1\leq
\varepsilon_1)$.
By our choice for $z'$, one has
\[
(n \geq N) \Rightarrow\biggl(\biggl\|\frac{y_n}{\|y_n\|
_1}-\frac
{z'}{\|z'\|_1}\biggr\|_1 \leq2\varepsilon_1\biggr)
\]
and, consequently,
$ \|y_n-h_n.z'\|_1 \leq2\varepsilon_1\|y_n\|_1 +\|z'\|_1$.
Thus, increasing $N$ if necessary, one has, for every $n\ge N$,
$ \|y_n-h_n.z'\|_1 \leq3\varepsilon_1\|y_n\|_1$.
But if $N$ is large enough, Lemma \ref{reglo} ensures that
\[
\forall n\ge N\qquad |\sigma(y_n)-\sigma(h_n.z')| \leq3C \varepsilon
_1\|
y_n\|_1.
\]
Finally, for every large $n$, we have
\begin{eqnarray*}
&&
|\sigma(y_n)-\mu(y_n)| \\
&&\qquad \leq |\sigma(y_n)-\sigma
(h_n.z')|+|\sigma
(h_n.z')-\mu(h_n.z')|
+ |\mu(h_n.z')-\mu(y_n)| \\
&&\qquad \leq 3C \varepsilon_1\|y_n\|_1 + h_n \biggl|\frac{\sigma
(h_n.z')}{h_n}-\mu(z')\biggr| + L\|h_n.z'-y_n\|_1 \\
&&\qquad \leq 3C \varepsilon_1\|y_n\|_1 +(1+\varepsilon_1)\frac{\|y_n\|
_1}{\|
z'\|_1}\biggl|\frac{\sigma(h_n.z')}{h_n}-\mu(z')\biggr|+
3\varepsilon_1
L\|y_n\|_1.
\end{eqnarray*}
But the a.s. convergence in the $z'$ direction ensures that for every
large $n$,
\[
\biggl|\frac{\sigma(h_n.z')}{h_n}-\mu(z')\biggr| \leq\varepsilon_1.
\]
Now if $\epsilon_1>0$ is small, we obtain, for every large $n$,
$|\sigma(y_n)-\mu(y_n)| < \varepsilon\|y_n\|_1;$
this brings contradiction and the proof is complete.\vspace*{-3pt}
\end{pf}

We can now prove the shape result for the ``fattened'' version $\tilde
G_t$ of $G_t=\{x\in\Zd\dvtx\sigma(x)\le t\}$; we recall that $A_\mu$ is
the unit ball for $\mu$.\vspace*{-3pt}
%
%
\begin{lemme}
\label{formesigma}
For every $\epsilon>0$, $\Pbarre$-a.s., for every large $t$,
\[
(1-\epsilon)A_{\mu}\subset\frac{\tilde G_t}t\subset(1+\epsilon
)A_{\mu}.\vspace*{-3pt}
\]
\end{lemme}
\begin{pf}
Let us prove by contradiction that if $t$ is large enough, $\frac
{G_t}t\subset(1+\epsilon)A_{\mu}$.
Thus assume that there exists an increasing sequence$(t_n)_{n\ge1}$,
with $t_n\to+\infty$ and $\frac{G_{t_n}}{t_n}\not\subset
(1+\epsilon
)A_{\mu}$;
so there exists $x_n$ with $\sigma(x_n)\le t_n$ and $\mu
(x_n)/t_n>1+\epsilon$.
So $\mu(x_n)/\sigma(x_n)>1+\epsilon$, which contradicts the uniform
convergence of Lemma \ref{convunifsigma}. Since $\mu
(x_n)>t_n(1+\epsilon
)$, the sequence $(\|x_n\|)_{n\ge1}$ goes to infinity.

For the\vspace*{-1pt} inverse inclusion, we still assume by contradiction that there
exists an increasing sequence $(t_n)_{n\ge1}$, with $t_n\to+\infty$
and $(1-\epsilon)A_{\mu}\not\subset\frac{\tilde G_{t_n}}{t_n}$; this
means we can find $x_n$ with $\mu(x_n)\le(1-\epsilon)t_n$, but
$\sigma
(x_n)>t_n$.
Since $t_n$ goes to $+\infty$, the sequence $(x_n)_{n\ge1}$ is not
bounded and satisfies $\frac{\mu(x_n)}{\sigma(x_n)}<1-\epsilon$; this
contradicts once again the uniform convergence of Lemma \ref
{convunifsigma} and the proof is complete.\vspace*{-3pt}
\end{pf}

Then we immediately recover the uniform convergence result for the
hitting time $t$ via Lemma \ref{pareil}, and, by an argument similar to
the one used in Lemma \ref{formesigma}, the asymptotic shape result for
the ``fattened'' version $\tilde H_t$ of $H_t=\{x\in\Zd\dvtx t(x)\le t\}$.\vspace*{-3pt}
%
%
\begin{lemme}
\label{etpourt}
$\Pbarre$-a.s., $ \lim_{\|x\|\to+\infty} \frac
{t(x)-\mu
(x)}{\|x\|}=0$,
and for every $\epsilon>0$, $\Pbarre$-a.s., for every large $t$,
$ (1-\epsilon)A_{\mu}\subset\frac{\tilde H_t}t\subset
(1+\epsilon)A_{\mu}$.\vspace*{-3pt}
\end{lemme}

It only remains now to prove the shape result for the coupled zone
$\tilde K'_t$, which is the ``fattened'' version of
$K'_t=\{x\in\Zd\dvtx
\forall s\ge t \xi^0_s(x)=\xi^{\Zd}_s(x)\}$.\vspace*{-3pt}
%
%
\begin{lemme}
\label{zoneC}
For every $\epsilon>0$, $\Pbarre$-a.s., for every large $t$,
$ (1-\epsilon)A_{\mu}\subset\frac{\tilde K'_t \cap
\tilde G_t}t$.\vspace*{-3pt}
\end{lemme}
\begin{pf}
Since $t\mapsto K'_t\cap G_t$ is nondecreasing, we use the same scheme
of proof as for Lemma \ref{formesigma}. We set, for $x \in\Zd$,
\[
t'(x)=\inf\{t\ge0\dvtx x\in K'_t\cap G_t\}.
\]
It is then sufficient to prove that $\Pbarre$-a.s.,
$\lim_{\|x\|\to+\infty} \frac{|t'(x)-\sigma(x)|}{\|x\|}=0$.\vadjust{\goodbreak}

By definition, $t'(x)\ge\sigma(x)$; thus it is sufficient to prove the
existence of constants $A',B'>0$ such that
%
%
\begin{equation}
\label{unbut}
\forall x\in\Zd,\forall s\ge0 \qquad\Pbarre\bigl( t'(x)-\sigma
(x)\ge
s\bigr)\le A'e^{-B' s}.
\end{equation}

$\bullet$ First note that for every $t \ge0$, $K_{\sigma
(x)+t}\supset
x+K_t \circ\tilde{\theta}_x$.

Indeed, let $z \in x+K_t \circ\tilde{\theta}_x$.
First consider the case $z \notin\xi_{\sigma(x)+t}^{\Zd}$.
Since, by additivity (\ref{additivite}), $\xi_{\sigma
(x)+t}^{0}\subset
\xi_{\sigma(x)+t}^{\Zd}$, we have $z \notin\xi_{\sigma(x)+t}^{0}$,
and so that $z\in K_{\sigma(x)+t}$.

Consider now the case $z\in\xi^{\Zd}_{\sigma(x)+t}$. Since, by
additivity, $\xi^{\Zd}_{\sigma(x)} \subset\xi^{\Zd}_{0}\circ
\tilde
{\theta}_x$, we have $y=z-x \in\xi^{\Zd}_{t}\circ\tilde{\theta}_x$.
But since $y\in K_t \circ\tilde{\theta}_x$, the definition of $K_t$
implies that $\xi^{0}_{t}(y)\circ\tilde{\theta}_x=\xi^{\Zd
}_{t}(y)\circ
\tilde{\theta}_x=1$.
Since\vspace*{-1pt} $x\in\xi^0_{\sigma(x)}$ and $y\in\xi^{0}_{t}\circ\tilde
{\theta
}_x$, we obtain $z=x+y\in\xi^0_{\sigma(x)+t}$, and so $z\in
K_{\sigma(x)+t}$.\vspace*{1pt}

$\bullet$ Fix $s \ge0$. The previous point says that
\[
\biggl(\bigcap_{t\ge s} K_{\sigma(x)+t}\biggr) \supset\biggl(
x+\bigcap
_{t\ge s}(K_t \circ\tilde{\theta}_x)\biggr)\quad \mbox{and so}\quad
K'_{\sigma(x)+s}\supset\bigl(x+(K'_s \circ\tilde{\theta}_x)\bigr).
\]
%
Since $\Pbarre$ is invariant under $\tilde{\theta}_x$, we get
\begin{eqnarray*}
\Pbarre\bigl(t'(x)>\sigma(x)+s\bigr) & = & \Pbarre\bigl(x\notin K'_{\sigma
(x)+s}\cap
G_{\sigma(x)+s}\bigr) \\
& =& \Pbarre\bigl(x\notin K'_{\sigma(x)+s}\bigr)
\le \Pbarre\bigl( x\notin(x+K'_s \circ\tilde{\theta}_{x})
\bigr) \\
& \le& \Pbarre(0\notin K'_s\circ\tilde{\theta}_{x} )=\Pbarre
(0\notin K'_s).
\end{eqnarray*}
We conclude with (\ref{petitsouscouple}).
\end{pf}

Let us finally prove that $\mu$ is a norm. Considering Lemma \ref
{munorme}, we only have to prove that $\mu(x)=\mu(-x)$ holds for each
$x\in\Zd$. This would be immediate if we had supposed that the law
$\nu
$ of the random environment was invariant under the central symmetry.
But, in general, we have to use a time reversal argument and the shape
theorem for the coupled zone.
We first give a characterization of $\mu$ that will allow us to use the
symmetries of the model.

%
\begin{lemme}
\label{reecriture}
Let us define $\check{\P}$ by $\check{\P}(A)=\int_{\Lambda} \P
_{\lambda
}(A) \,d\nu(\lambda)$. Then, for each \mbox{$x\in\Zd$}
\[
\mu(x)=\sup\Bigl\{a>0; \liminf_{n\to+\infty}\check{\P}(nx\in\xi
_{n/a})>0\Bigr\}.
\]
\end{lemme}
\begin{pf}
Define $g(x)=\sup\{a>0; \liminf_{n\to+\infty}\check{\P}(nx\in
\xi_{n/a})>0\}$.
Let $a>\mu(x)$. By the asymptotic shape theorem,
$\lim_{n\to+\infty}\P_{\lambda}(nx\in\xi_{n/a})=0$ $\nu$
almost surely holds, so by dominated convergence, $\lim_{n\to
+\infty}\check{\P}(nx\in\xi_{n/a})=0$. This gives $g(x)\le\mu(x)$.
Now take some $a<\mu(x)$. We will show that
\[
\liminf_{n\to+\infty}{\P}_{\lambda}(nx\in\xi_{n/a})\ge\P
_{\lambda_{\min}}(\tau=+\infty)^2,\qquad \nu\mbox{-a.s.},
\]
which will give
\[
\liminf_{n\to+\infty}\check{\P}(nx\in\xi_{n/a})\ge\int
_{\Lambda} \liminf_{n\to+\infty}{\P}_{\lambda}(nx\in\xi
_{n/a}) \,d\nu\ge\P_{\lambda_{\min}}(\tau=+\infty)^2
\]
by Fatou's lemma, whence $g(x)\ge a$, which will lead to $g(x)=\mu(x)$.
Obviously,
$\liminf_{n\to+\infty}{\P}_{\lambda}(nx\in\xi_{n/a})\ge
\liminf_{n\to+\infty}{\P}_{\lambda}(nx\in\xi_{n/a},\tau
^0=+\infty)$.
However, the convergence theorem for the coupled zone implies that
\[
\lim_{n\to+\infty}\P_{\lambda}(nx\in\xi^{\Zd}_{n/a},nx\notin
\xi_{n/a},\tau^0=+\infty)=0, \qquad\nu\mbox{-a.s.},
\]
hence, by a classical time-reversal argument and using the FKG inequality,
\begin{eqnarray*}
\liminf_{n\to+\infty}{\P}_{\lambda}(nx\in\xi_{n/a})&\ge
&\liminf_{n\to+\infty}\P_{\lambda}(nx\in\xi^{\Zd}_{n/a},\tau
^0=+\infty)\\
&\ge&\liminf_{n\to+\infty}\P_{\lambda}(nx\in\xi^{\Zd
}_{n/a})\P_{\lambda}(\tau^0=+\infty)\\
& \ge&\liminf_{n\to+\infty}\P_{\lambda}\biggl(\tau^{nx}\ge\frac
{n}{a}\biggr)\P_{\lambda}(\tau^0=+\infty)\\
& \ge& \P_{\lambda_{\min}}(\tau=+\infty)^2,
\end{eqnarray*}
which ends the proof of the lemma.
\end{pf}

We now have a handsome expression to prove the symmetry property.
Actually, for every $x\in\Zd,t>0,\lambda\in\Lambda$, a time-reversal
argument proves that $\P_{\lambda}(x\in\xi_t)=\P_{x.\lambda
}(-x\in\xi
_t)$, hence integrating with respect to $\nu$ and using the invariance
of $\nu$ under the translation by $x$,
\[
\check{\P}(x\in\xi_t)=\check{\P}(-x\in\xi_t),
\]
which, with Lemma \ref{reecriture}, gives the symmetry of $\mu$.

\section{Uniform controls of the growth}
\label{restartestimees}

The aim of this section is to establish some of the uniform controls
announced in Proposition \ref{propuniforme}.
To control the growth of the contact process, we need some lemmas on
the Richardson model.

\subsection{Some lemmas on the Richardson model}
We call Richardson model with parameter $\lambda$ the time-homogeneous,
$\mathcal{P}(\Zd)$-valued Markov process\break $(\eta_t)_{t\ge0}$, whose
evolution is defined as follows: an empty site $z$ becomes infected at
rate $ \lambda\sum_{\|z-z'\|_1=1} \eta_t(z') $, the
different evolutions being independent.
Thanks\vspace*{1pt} to the graphical construction, we can, for each $\lambda\in
\Lambda$,
build a coupling of the contact process in environment $\lambda$
with the Richardson model with parameter $\lambda_{\max}$, in the
following way: at any time $t$, the space occupied by the contact
process is contained in the space occupied by the Richardson
model.\vadjust{\goodbreak}

The first lemma, whose proof is omitted, easily follows from the
representation of the Richardson model in terms of first passage
percolation, together with a path counting argument.\vspace*{-3pt}
%
%
\begin{lemme}
\label{Richardson1}
For every $\lambda>0$, there exist constants $A,B>0$ such that
\[
\forall t\ge0\qquad \P(\eta_1\not\subset B_t)\le A\exp(-B t).\vspace*{-3pt}
\]
\end{lemme}
%
%
\begin{lemme}
\label{Richardson2}
For every $\lambda>0$, there exist constants $A,B,M>0$ such that
\[
\forall s\ge0 \qquad\P(\exists t\ge0\dvtx\eta_t \not\subset
B_{Mt+s})\le A\exp(-Bs).\vspace*{-3pt}
\]
\end{lemme}
\begin{pf}
The representation of the Richardson model in terms of first passage
percolation ensures the existence of $A',B',M'>0$ such that for each $t
\ge0$,
%
%
\begin{equation}
\label{voirkesten}
\P(\eta_t\not\subset B_{M't}) \le A'\exp(-B't).
\end{equation}
For more details, one can refer to Kesten \cite{kesten}.

We first control the process in integer times thanks to the following estimate:
%
%
\begin{eqnarray}\label{olala}
\P(\exists k\in\N\dvtx\eta_k\not\subset B_{M'k+s/2})
& \le&\P\bigl(\exists k\in\N\dvtx\eta_{k+s/(2M')}\not\subset B_{M'k+s/2}\bigr)
\nonumber\\
&{\le} & \sum_{k=0}^{+\infty}\P\bigl( \eta_{k+s/(2M')}\not\subset
B_{M'k+s/2}\bigr) \\
& \le& \frac{A'}{1-\exp(-B')}\exp\biggl(-\frac{B's}{2M'}\biggr).
\nonumber
\end{eqnarray}
Let us now control the fluctuations between integer times.
Let $M>M'$,
%
%
\begin{eqnarray}\label{olili}
&& \P(\{\exists t\ge0\dvtx\eta_t\not\subset B_{Mt+s}\}\cap\{\forall k
\in\N, \eta_k\subset B_{M'k+s/2} \}) \nonumber\\[-8pt]\\[-8pt]
&&\qquad \le \sum_{k=0}^{+\infty}\P( \exists t\in[k,k+1]\dvtx\eta
_{k}\subset B_{M'k+s/2} \mbox{ and } \eta_t\not\subset B_{Mt+s}).
\nonumber
\end{eqnarray}
Then, denoting by $C'>0$ a constant such that $\Card{B_t}\le C'(1+t)^d$
and by $A,B$ the constants appearing in Lemma \ref{Richardson1},
%
%
\begin{eqnarray}\label{rouge}
&& \P( \exists t\in[k,k+1]\dvtx\eta_{k}\subset B_{M'k+s/2} \mbox{ and
} \eta_t\not\subset B_{Mt+s})\nonumber\\
&&\qquad \le \P(\eta_{k}\subset B_{M'k+s/2} \mbox{ and } \eta_{k+1}\not
\subset B_{Mk+s})\nonumber\\
&&\qquad \le \Card{B_{M'k+s/2}} \P\bigl(\eta_1\not\subset B_{ k(M-M')+s/2}\bigr)
\\
&&\qquad \le C'(1+ M'k+s/2)^dA\exp\bigl(-B\bigl(k(M-M')+s/2\bigr)\bigr)\nonumber\\
&&\qquad \le AC' (1+s/2)^d\exp(-Bs/2)(1+M'k)^d\exp\bigl(-B\bigl(k(M-M')\bigr)\bigr)\nonumber.
\end{eqnarray}
Inequality (\ref{rouge}) comes from the Markov property and from the
subadditivity of the contact process.
Since the series $((1+M'k)^d\exp(-B(k(M-M')))_{k\ge1}$ converges, the
desired result follows from (\ref{olala}) and (\ref{olili}).\vadjust{\goodbreak}
\end{pf}

\subsection{A restart procedure}

We will use here a so-called restart argument, which can be summed up
as follows. We couple the system that we want to study (the strong
system) with a system that it stochastically dominates (the weak
system), and that is best understood.
Then, we can transport some of the properties of the known system to
the one we study; we let the processes simultaneously evolve and, each
time the weaker dies and the stronger remains alive, we restart a copy
of the weakest, coupled with the strongest again.
Thus, either both processes die before we found any weak process
surviving. In this case, the control of large finite lifetimes for the
weak can be transposed to the strongest one, or the strongest
indefinitely survives and is finally coupled with a weak surviving one.
In that case, a~bound for the time that is necessary to find a
successful restart permits us to transfer properties of the weak
surviving process to the strong one.

This technique is already old; that can be found, for example, in
Durrett~\cite{MR757768}, Section 12, in a very pure form. It is also
used by Durrett and Griffeath \cite{MR656515}, in order to transfer
some controls for the one-dimensional contact process to the contact
process in a larger dimension. We will use it here by coupling the
contact process in inhomogeneous environment $\lambda\in\Lambda$ with
the contact process with a constant birth rate $\lambda_{\min}$. Here,
the assumption $\lambda_{\min}>\lambda_c(\Zd)$ matters.

To this end, we will couple collections of Poisson point processes. Fix
$\lambda\in\Lambda$. We can build a probability measure $\tilde\P
_{\lambda}$ on $\Omega\times\Omega$ under which:
\begin{itemize}
\item the first coordinate $\omega$ is a collection $((\omega_e)_{e
\in
\Ed}, (\omega_z)_{z \in\Zd})$ of Poisson point processes, with
respective intensities $(\lambda_e)_{e \in\Ed}$ for the bond-indexed
processes, and intensity $1$ for the site-indexed processes.
\item the second coordinate $\eta$ is a collection $((\eta_e)_{e \in
\Ed
}, (\eta_z)_{z \in\Zd})$ of Poisson point processes, with intensity
$\lambda_{\min}$ for the bond-indexed processes, and intensity $1$ for
the site-indexed processes.
\item site-indexed Poisson point processed (death times) coincide; for
every \mbox{$z \in\Zd$}, $\eta_z=\omega_z$.
\item bond-indexed Poisson point processed (birth-times candidates) are
coupled; for each $e \in\Ed$, the support of $\eta_e$ is included in
the support of $\omega_e$.
\end{itemize}
We denote by $\xi^A=\xi^A(\omega,\eta)$ the contact process in environment
$\lambda$ starting from $A$ and built from the Poisson process
collection $\omega$, and
$\zeta^B=\zeta^B(\omega,\eta)$ the contact process in environment
$\lambda_{\min}$ starting from $B$ and built from the Poisson process
collection $\eta$.
If $B \subset A$, then $\tilde\P_{\lambda}$ almost surely, $\zeta_t^B
\subset\xi_t^A$
holds for each $t \ge0$. We can note that the process $(\xi^A,\zeta
^B)$ is a Markov process.

We introduce the lifetimes of both processes:
\[
\tau= \inf\{t\ge0\dvtx\xi^0_t=\varnothing\} \quad\mbox{and}\quad\mbox{for }x\in
\Zd\qquad
\tau'_x = \inf\{t\ge0\dvtx\zeta^x_t=\varnothing\}.
\]
Note that the law of $\tau'_x$ under $\tilde{\P}_{\lambda}$ is the law
of $\tau_x$ under $\P_{\lambda_{\min}}$; it actually does not
depend on
the process starting point, because the model with constant birth rate
is translation invariant.

We recursively define a sequence of stopping times $(u_k)_{k \ge0}$
and a sequence of points $(z_k)_{k \ge0}$, letting $u_0=0$, $z_0=0$
and for each $k \ge0$:
\begin{itemize}
\item if $u_k<+\infty$ and $\xi_{u_k}\neq\varnothing$, then
$u_{k+1}=\tau'_{z_k} \circ\theta_{u_k}$;
\item if $u_k=+\infty$ or if $\xi_{u_k}= \varnothing$, then
$u_{k+1}=+\infty$;
\item if $u_{k+1}<+\infty$ and $\xi_{u_{k+1}}\neq\varnothing$, then
$z_{k+1}$ is the smallest point of $\xi_{u_{k+1}}$ for the
lexicographic order;
\item if $u_{k+1}=+\infty$ or if $\xi_{u_{k+1}}= \varnothing$, then
$z_{k+1}=+\infty$.
\end{itemize}
In other words, until $u_k<+\infty$ and $\xi_{u_k}\neq\varnothing$, we
take in
$\xi_{u_k}$ the smallest point $z_k$ for the lexicographic order, and
look at the lifetime of the weakest process, namely, $\zeta$, starting
from $z_k$ at time $u_k$. The restart procedure can stop in two ways;
either we find $k$ such that $u_k<+\infty$ and $\xi_{u_k}=
\varnothing
$, which implies that the strongest process (which contains the weak)
precisely dies at time $u_k$,
or we find $k$ such that $u_k<+\infty$, $\xi_{u_k} \neq\varnothing$
and $u_{k+1}=+\infty$. In this case, we have found a point $z_k$ such
that the weak process which starts from $z_k$ at time $u_k$ survives;
particularly, this implies that the strongest also survives.
We then define
\[
K=\inf\{n\ge0\dvtx u_{n+1}=+\infty\}.
\]
The name of the $K$ variable is chosen by analogy with Section \ref{sigma}.
The current section being independent from the rest of the article,
confusion should not be possible. It comes from the preceding
discussion that
%
%
\begin{equation}
\label{arev}
(\tau=+\infty\Longleftrightarrow\xi^0_{u_K}\neq\varnothing)
\quad\mbox{and if}\quad \tau<+\infty\qquad \mbox{then } u_K=\tau.
\end{equation}

We regroup in the next lemma some estimates on the restart procedure
that are necessary to prove Proposition \ref{propuniforme}. Recall that
$\rho$ is introduced in (\ref{uniftau}).

%
\begin{lemme}
\label{restart}
We work in the preceding frame. Then:
\begin{itemize}
\item$\forall\lambda\in\Lambda,\forall n\in\N$ $\tilde
{\P
}_\lambda(K>n)\le(1-\rho)^n$.
\item$\forall B \in\mathcal B(\mathcal D)$ $\tilde{\P}_\lambda
(\tau
=+\infty, \zeta^{z_K}\circ\theta_{u_K} \in B) = \P_\lambda(\tau
=+\infty
)\Pbarre_{\lambda_{\min}}(\xi^0 \in B)$.
\item there exist $\alpha,\beta>0$ such that for every $\lambda\in
\Lambda$, $\tilde{\E}_{\lambda}(\exp(\alpha u_K))<\beta$.
\end{itemize}
\end{lemme}
\begin{pf}
By the strong Markov property, we have
\begin{eqnarray*}
\tilde{\P}_\lambda(K\ge n+1) & = & \tilde{\P}_\lambda
(u_{n+1}<+\infty)
\\
& = & \tilde{\P}_\lambda(u_n<+\infty, \xi_{u_n} \neq
\varnothing,
\tau'_{z_n} \circ\theta_{u_n}<+\infty) \\
& \le& \tilde{\P}_{\lambda}(u_n<+\infty) (1-\rho)\\
&=&\tilde{\P
}_\lambda(K
\ge n)(1-\rho).
\end{eqnarray*}
Thus, $K$ has a subexponential tail, which proves the first point.
Particularly,~$K$ is almost surely finite.

Using (\ref{arev}) and the strong Markov property, we also have
\begin{eqnarray*}
&& \tilde{\P}_\lambda(\tau=+\infty, \zeta^{z_K}\circ{\theta
}_{u_K} \in B)\\
&&\qquad= \tilde{\P}_\lambda(\xi_{u_K}\neq\varnothing, \zeta^{z_K}\circ
{\theta}_{u_K} \in B) \\
&&\qquad = \sum_{k=0}^{+\infty}\sum_{z \in\Zd} \tilde{\P}_\lambda
(K=k,
\xi^0_{u_k}\neq\varnothing, z_k=z, \zeta^{z_K}\circ{\theta
}_{u_K} \in B) \\
&&\qquad = \sum_{k=0}^{+\infty}\sum_{z \in\Zd} \tilde{\P}_\lambda
(u_k<+\infty, \xi^0_{u_k}\neq\varnothing, z_k=z, \tau'_{z_K}
\circ\theta_{u_k}=+\infty, \zeta^{z_K}\circ{\theta}_{u_K} \in
B) \\
&&\qquad = \sum_{k=0}^{+\infty}\sum_{z \in\Zd} \tilde{\P}_\lambda
(u_k<+\infty, \xi^0_{u_k}\neq\varnothing, z_k=z)\P_{\lambda
_{\min
}}(\tau=+\infty, \xi^{0} \in B) \\
&&\qquad = \P_{\lambda_{\min}}( \tau=+\infty, \xi^{0} \in B) \sum
_{k=0}^{+\infty}\tilde{\P}_\lambda(u_k<+\infty, \xi
^0_{u_K}\neq
\varnothing).
\end{eqnarray*}
Taking for $B$ the whole set of trajectories, we can identify
\[
\tilde{\P} (\tau=+\infty)= \P_\lambda(\tau=+\infty)=\P
_{\lambda_{\min
}}(\tau=+\infty)\sum_{k=0}^{+\infty}\tilde{\P}_\lambda
(u_k<+\infty,\xi
^0_{u_K}\neq\varnothing),
\]
which gives us the second point.

Since $\lambda_{\min}>\lambda_c(\Zd)$, the results by Durrett and
Griffeath \cite{MR656515} for large~$\lambda$, extended to the whole
supercritical regime by Bezuidenhout and Grimmett~\cite{MR1071804},
ensure the existence of $A,B>0$ such that
\[
\forall t\ge0\qquad \P_{\lambda_{\min}}(t \le\tau<+\infty) \le A
\exp(-Bt),
\]
which gives the existence of exponential moments for $\tau\one_{\{
\tau
<+\infty\}}$.
Since\break $\P_{\lambda_{\min}}(\tau=+\infty)>0$, we can choose (e.g., by
dominated convergence) some $\alpha>0$ such that $\E_{\lambda_{\min
}}(\exp(\alpha\tau)\one_{\{\tau<+\infty\}})=r<1$.

For $k\ge0$, we note
\[
S_k=\exp\Biggl( \alpha\sum_{i=0}^{k-1}\tau'_{z_i} \circ\theta_{u_i}
\Biggr)\one_{\{u_{k}<+\infty\}}.
\]
We note that $S_k$ is $\mathcal{F}_{u_{k}}$-measurable. Let $k\ge0$.
We have
\[
\exp(\alpha u_K)\one_{\{K= k\}}\le S_{k}.
\]
Thus, applying the strong Markov property at time $u_{k-1}<+\infty$, we
get, for $k \ge1$,
\begin{eqnarray*}
\tilde{\E}_\lambda\bigl[\exp(\alpha u_K)\one_{\{ K= k\}}\bigr] & \le&
\tilde{\E
}_\lambda(S_k)=\tilde{\E}_\lambda(S_{k-1}){\E}_{\lambda_{\min
}}\bigl(\exp
(\alpha\tau)\one_{\{\tau<+\infty\}}\bigr) \\
& \le& r\tilde{\E}_\lambda(S_{k-1}).
\end{eqnarray*}
Since $r<1$, it comes that $\tilde{\E}_\lambda[\exp(\alpha u_K)]\le
\frac{r}{1-r}<+\infty$.\vadjust{\goodbreak}
\end{pf}

\subsection{\texorpdfstring{Proof of Proposition \protect\ref{propuniforme}}{Proof of Proposition 5}}

Estimates (\ref{richard}) and (\ref{uniftau}) follow from a simple
stochastic comparison.
\begin{pf*}{Proof of (\ref{uniftau})}
It suffices to note that for every environment $\lambda\in\Lambda$
and each $z\in\Zd$, we have
\[
\P_{\lambda}(\tau^z=+\infty) \ge\P_{\lambda_{\min}}(\tau
^z=+\infty)=\P
_{\lambda_{\min}}(\tau^0=+\infty)>0.
\]
\upqed\end{pf*}
\begin{pf*}{Proof of (\ref{richard})}
We use the stochastic domination of the contact process in environment
$\lambda$ by the Richardson model with parameter $\lambda_{\max}$. For
this model, (\ref{voirkesten})~ensures a growth which is at least
linear.
\end{pf*}

Then, it remains to prove (\ref{grosamasfinis}), (\ref{retouche}) and
(\ref{petitsouscouple}) with a restart procedure.
\begin{pf*}{Proof of (\ref{grosamasfinis})}
Let $\alpha,\beta>0$ as given in the third point of Lemma \ref
{restart}. Recall that $u_{K}=\tau$ on $\{\tau<+\infty\}$.
For each $\lambda\in\Lambda$ and each $t>0$, we have
\begin{eqnarray*}
\P_\lambda(t<\tau<+\infty)& =& \P_\lambda(e^{\alpha t}<e^{\alpha
\tau
},\tau<+\infty)
= \tilde{\P}_\lambda(e^{\alpha t}<e^{\alpha u_{K}},\tau<+\infty)\\
& \le& \tilde{\P}_\lambda(e^{\alpha t}<e^{\alpha u_{K} })\le
e^{-\alpha
t} \tilde{\E}_\lambda e^{\alpha u_{K}} \le\beta e^{-\alpha t},
\end{eqnarray*}
which concludes the proof.
\end{pf*}
\begin{pf*}{Proof of (\ref{retouche})}
Since $\lambda_{\min}>\lambda_c(\Zd)$, Durrett and Griffeath's
results~\cite{MR656515} for large $\lambda$, extended to the whole
supercritical regime by Bezuidenhout and Grimmett~\cite{MR1071804},
ensure the existence of constants $A,B,c>0$ such that, for each $y \in
\Zd$, for each $t \ge0$,
%
%
\begin{equation}
\label{blibli}
\Pbarre_{\lambda_{\min}}\biggl( t(y) \ge\frac{\|y\|}{c}+t \biggr)
\le A
\exp(-Bt).
\end{equation}
Besides, the domination by the Richardson model with parameter $\lambda
_{\max}$ and Lem\-ma~\ref{Richardson2} ensure the existence of $A,B,M>0$
such that for every $\lambda\in\Lambda$, for each $s \ge0$,
%
%
\begin{equation}
\label{blabla}
\P_{\lambda}(\exists t\ge0, \xi^0_t \not\subset B_{Mt+s}) \le
A\exp(-Bs).
\end{equation}
By decreasing $c$ or increasing $M$ if necessary, we can also assume
that $\frac{c}{M}\le1$. Now,
\begin{eqnarray*}
&& \tilde{\P}_{\lambda} \biggl(t(y) \ge\frac{\|y\|}c +t, \tau
=+\infty
\biggr) \\
&&\qquad \le \tilde{\P}_{\lambda} \biggl( u_K \ge\frac{tc}{6M} \biggr)+
\tilde{\P}_{\lambda} \biggl( u_K \le\frac{tc}{6M}, \xi^0_{u_K}
\not
\subset B_{tc/3} \biggr) \\
&&\qquad\quad{} +\tilde{\P}_{\lambda} \biggl(\tau=+\infty, u_K \le\frac{tc}{6M},
\xi^0_{u_K} \subset B_{tc/3}, t(y) \ge\frac{\|y\|}c +t\biggr).
\end{eqnarray*}
By Lemma \ref{restart}, $u_K$ has exponential moments, so we can bound
the first term; there exist $C,\alpha>0$ such that for each $\lambda
\in\Lambda$, for each $t\ge0$,
\[
\tilde{\P}_{\lambda} \biggl(u_K \ge\frac{tc}{6M}\biggr) \le C\exp
\biggl(-\frac{\alpha c t}{6M}\biggr).
\]
The second term is controlled with the help of (\ref{blabla}):
\[
\tilde{\P}_{\lambda} \biggl(u_K \le\frac{tc}{6M}, \xi^0_{u_K}
\not
\subset B_{tc/3} \biggr) \le\P_{\lambda}\bigl(\exists t\ge0, \xi^0_t
\not\subset B_{Mt+({tc})/6}\bigr) \le A\exp\biggl(-B\frac{tc}6\biggr).
\]
It remains to bound the last term. We note here
\[
t'(y)=\inf\{t \ge0\dvtx y \in\zeta^0_t\}.
\]
Recall that if $\tau=+\infty$, then $\xi_{u_K} \neq\varnothing$ and
$z_K$ is well defined.
Since $t(y)$ is the hitting time of $y$ and $\xi^0_t\supset\zeta^0_t $
for each $t$, we have, on $\{\tau=+\infty\}$,
\[
t(y)\le u_{K}+t'(y-z_K) \circ T_{z_K} \circ\theta_{u_K}.
\]
If $u_K \le\frac{tc}{6M} \le\frac{t}6$, then
$t(y) \le\frac{t}{6}+t'(y-z_K)\circ T_{z_K} \circ\theta_{u_K}$.
If, moreover, $\xi^0_{u_K} \subset B_{tc/3}$, we have $\|y \| \ge\|
y-z_K\|-\frac{tc}3$,
which gives, with the second point in Lemma \ref{restart},
\begin{eqnarray*}
&& \tilde{\P}_{\lambda} \biggl(\tau=+\infty, u_K \le\frac{ct}{6M},
\xi^0_{u_K} \subset B_{ct/3}, t(y) \ge\frac{\|y\|}c +t\biggr)
\\
&&\qquad \le \tilde{\P}_{\lambda}\biggl(\tau=+\infty, 
t'(y-z_K) \circ T_{z_K} \circ\theta_{u_K} \ge\frac{\|y-z_K\|}c
+\frac{t}{2}\biggr) \\
&&\qquad \le \P_\lambda(\tau=+\infty) \sup_{z \in\Zd}\Pbarre
_{\lambda_{\min
}}\biggl(t(y-z) \ge\frac{\|y-z\|}c +\frac{t}2\biggr)
\le A\exp(-Bt/2),
\end{eqnarray*}
where the last inequality follows from (\ref{blibli}). The proof is complete.
\end{pf*}
\begin{pf*}{Proof of (\ref{petitsouscouple})}
Let $s \ge0$, and denote by $n$ the integer part of $s$. Let $\gamma
>0$ be a fixed number, whose precise value will be specified later:
\begin{eqnarray*}
\Pbarre(0\notin K'_s)
&=& \Pbarre(\exists t \ge s\dvtx0 \notin K_t) \\
& \le& \sum_{k=n}^{+\infty}\Pbarre(B_{\gamma k}\not\subset K_k)\\
&&{}+
\sum_{k=n}^{+\infty}\Pbarre\bigl(B_{\gamma k}\subset K_k, \exists t\in
[k,k+1) \mbox{ such that } 0\notin K_t\bigr).
\end{eqnarray*}

Let us first bound the second sum. Fix $k \ge n$. Assume that
$B_{\gamma k}\subset K_k$ and consider $t\in[k,k+1) \mbox{ such that }
0\notin K_t$. Then, there exists $x\in\Zd$ such that $0\in\xi
^x_t\setminus\xi^0_t$. Since $0\in\xi^x_t$ and $t \ge k$, there
exists $y \in\Zd$ such that $y \in\xi^x_k$ and $0 \in\xi^y_{t-k}
\circ\theta_k$.\vadjust{\goodbreak} If $y \in B_{\gamma k}\subset K_k$, then $\xi
^{0}_k(y)=\xi^{\Zd}_k(y)=1$, which implies that $y \in\xi^{0}_k$. Now,
since $0 \in\xi^y_{t-k} \circ\theta_k$, we obtain
$0\in\xi^0_t$, which contradicts the assumption $0\notin\xi^0_t$.
Thus, we necessarily have $y \notin B_{\gamma k}$, so
%
\begin{eqnarray*}
&& \Pbarre_{\lambda}\bigl(B_{\gamma k}\subset K_k, \exists t\in
[k,k+1)\mbox
{ tel que } 0\notin K_t\bigr) \\
&&\qquad \le \frac1{\P_\lambda(\tau=+\infty)}\P_{\lambda}
\biggl(\theta_k^{-1}
\biggl( 0\in\bigcup_{s\in[0,1]} \xi^{\Zd\setminus B_{\gamma
k}}_s \biggr) \biggr) \\
&&\qquad \le \frac1\rho\P_{\lambda}\biggl(0\in\bigcup_{s\in[0,1]} \xi
^{\Zd\setminus B_{\gamma k}}_s\biggr)
\le \frac1\rho\P_{\lambda}\biggl(\bigcup_{s\in[0,1]} \xi
^{0}_s\not\subset B_{\gamma k}\biggr)\\
&&\qquad = \frac1{\rho}\P_{\lambda
}(H_1^0\not\subset B_{\gamma k}).
\end{eqnarray*}
Since the Richardson model with parameter $\lambda_{\max}$
stochastically dominates the contact process in environment $\lambda$,
we control the last term thanks to Lemma \ref{Richardson1}.

To control the first sum, it is sufficient to prove that there exist
positive constants $A,B,\gamma$ (and this will fix the precise value of
$\gamma$) such that for each
$\lambda\in\Lambda$ and each $ t\ge0$
%
%
\begin{equation}\label{souscouple}
\P_{\lambda}(B_{\gamma t}\not\subset K_t, \tau^0=+\infty) \le
A\exp(-B t).
\end{equation}
%

The number of integer points in a ball being polynomial with respect to
the radius, it is sufficient to prove that there exist some constants
$A,B,c'>0$ such that for each $t\ge0$, for each $x \in\Zd$,
%
%
\begin{equation}
\label{monbutla}
\|x\|\le c't \quad\Longrightarrow\quad\tilde{\P}_{\lambda}(\xi^{0}_t\ne
\varnothing,x\in\xi^{\Zd}_t\setminus\xi^{0}_t)\le A\exp(-Bt).
\end{equation}
To prove (\ref{monbutla}), we will use the following result, that has
been obtained by Durrett \cite{MR1117232} as a consequence of the
Bezuidenhout and Grimmett construction \cite{MR1071804}. If $\xi^0$ and
$\tilde{\xi}^x$ are two independent contact processes with parameter
$\lambda>\lambda_c(\Zd)$, respectively, starting from $0$ and from
$x$, then there exist positive constants $A,B,\alpha$ such that for
each $t\ge0$ and each $x \in\Zd$,
%
%
\begin{equation}
\label{durrett}
\|x\|\le\alpha t \quad\Longrightarrow\quad\P(\xi^0_t\cap\tilde{\xi
}^x_t=\varnothing, \tilde{\xi}^x_t\ne\varnothing,{\xi}^0_t\ne
\varnothing)\le A\exp(-Bt).
\end{equation}
Let $\alpha$ and $M$ be the constants, respectively, given by
equations (\ref{durrett}) and~(\ref{richard}).
We put $c'=\alpha/2$ and choose $\epsilon>0$ such that $c'+2\epsilon
M\le\alpha$.

Let $a\in B^0_{\alpha t/4}$ and $b\in B^x_{\alpha t/4}$. We set
\[
\alpha_{a,s}=\zeta^a_s\circ\theta_{\epsilon t/2}\quad\mbox{and}\quad\beta
_{b,s}=\bigl\{y\in\Zd\dvtx b\in\zeta^{y}_s\circ\theta_{t(1-\epsilon
/2)-s}\bigr\}.
\]
Then, $(\alpha_{a,s})_{0\le s\le t/2(1-\epsilon)}$ and $(\beta
_{a,s})_{0\le s\le t/2(1-\epsilon)}$
are independent contact processes with constant birth rate $\lambda
_{\min}$, respectively, starting from $a$ and from~$b$.
The process $(\beta_{a,s})_{0\le s\le t/2(1-\epsilon)}$ is a contact
process, but for which the time axis has been reverted.
In the same way, we set
\[
\hat{\xi}_s^x=\{y\in\Zd\dvtx x\in\xi_s^y\circ\theta_{t-s}\}.\vadjust{\goodbreak}
\]
Note that $(\hat{\xi}^x_{s})_{0\le s\le t/2}$ has the same law as
$({\xi}^x_{s})_{0\le s\le t/2}$.
Note that:
\begin{itemize}
\item assuming $a\in\xi^0_{\epsilon t/2}$, $\alpha_{a,(1-\epsilon)
t/2}\cap\beta_{b,(1-\epsilon)t/2}\ne\varnothing$
and
$b\in\hat{\xi}^x_{\epsilon t/2}$, then $x\in\xi_t^0$;
\item if $x\in\xi_t^{\Zd}$, then $\hat{\xi}^x_{t/2}$ is nonempty;
\item if ${\xi}^0_{t}$ is nonempty, then ${\xi}^0_{t/2}$ is nonempty.
\end{itemize}
Thus, letting
\[
E^0  =  \{\xi^0_{t/2}\ne\varnothing\}\setminus
\bigl\{
\exists a\in B^0_{\alpha t/4} \cap\xi^0_{\epsilon t/2}\dvtx
\alpha_{a,(1-\epsilon) t/2}\ne\varnothing
\bigr\}
\]
and
\[
\hat{E}^x = \{\hat{\xi}^x_{t/2}\ne
\varnothing
\}\setminus
\bigl\{
\exists b\in B^x_{\alpha t/4} \cap\hat{\xi}^x_{\epsilon t/2}\dvtx
\beta_{b,(1-\epsilon) t/2}\ne\varnothing
\bigr\},
\]
we get
%
%
\begin{eqnarray}
\label{decouplage}
\tilde{\P}_{\lambda}(\xi^{0}_t\ne\varnothing, x\in\xi^{\Zd
}_t\setminus\xi^{0}_t)
& \le& \tilde{\P}_{\lambda}( \xi^{0}_{t/2} \ne\varnothing,
\hat
{\xi}^x_{t/2} \ne\varnothing, \xi^{0}_{t/2}\cap\hat{\xi
}^x_{t/2}=\varnothing) \nonumber\\[-8pt]\\[-8pt]
& \le& \tilde{\P}_{\lambda}(E^0)+\tilde{\P}_{\lambda}(\hat
{E}^x) +S
\nonumber,
\end{eqnarray}
where\vspace*{-1pt}
$S
= \sum_{a\in B^0_{\alpha t/4},b\in B^x_{\alpha t/4}}\tilde
{\P}_{\lambda}(\alpha_{a,{(1-\epsilon) t}/2}\ne
\varnothing,
\beta_{b,{(1-\epsilon) t}/2}\ne\varnothing, \alpha_{a,
{(1-\epsilon) t}/2}\cap\beta_{b,{(1-\epsilon) t}/2}=\varnothing
)$.\vspace*{1pt}

For every couple $(a,b)$ that appears in $S$, we have
$\|a-b\|\le\|a\|+\|b-x\|+\|x\|\le\alpha t/4+ \alpha t/4+\alpha
t/2=\alpha t$, which allows us to use (\ref{durrett}), and gives the
existence of constants $A,B,C'>0$ such that
\[
S \le C'(1+\alpha t/4)^{2d}A\exp\bigl(-B(1-\epsilon)t/2\bigr).
\]
By another time reversal, we see that $\tilde{\P}_{\lambda}(\hat
{E}^x)=\tilde{\P}_{\T{x}{\lambda}}(E^0)$; then it suffices to control
$\tilde{\P}_{\lambda}(E^0)$ uniformly in $\lambda$.
Let
\[
E_1=\{\xi^{0}_{t/2}\ne\varnothing\}\setminus\bigl\{\exists a\in\Zd\dvtx
a\in
\xi^0_{\epsilon t/2}, \alpha_{a,(1-\epsilon) t/2}\ne\varnothing\bigr\}.
\]
We have
$\tilde{\P}_{\lambda}(E^0)\le\tilde{\P}_{\lambda}(E_1)+\tilde
{\P
}_{\lambda}(\xi^0_{\epsilon t/2}\not\subset B^0_{\alpha t/4})$.
By the choice we made for $\epsilon$ and inequality (\ref{richard}),
we have
\[
\forall\lambda\in\Lambda,\forall t\ge0 \qquad\tilde{\P
}_{\lambda
}\bigl(\xi^0_{\epsilon t/2}\not\subset B(0,\alpha t/4)\bigr)\le A\exp
(-B\epsilon t/2).
\]
Thanks to the restart Lemma \ref{restart}, we can see that
\[
\tilde{\P}_{\lambda}(u_K> \epsilon t/2)\le\beta\exp(-\alpha
\epsilon t/2).
\]
Suppose then that $u_K\le\epsilon t/2$ and $\xi^{0}_{t/2}\ne
\varnothing
$: $z_K$ is thus well defined and we have
$\tau'_{z_K} \circ\theta_{u_K}=+\infty$.
Then, there exists an infinite infection branch in the coupled process
in environment $\lambda_{\min}$ starting from $\xi_{u_K}^0$. This branch
contains at least one point $a\in\xi^0_{(1-\epsilon)t/2}$.
By construction $a\in\xi^0_{(1-\epsilon)t/2}$ and
$\alpha_{a,(1-\epsilon) t/2}\ne\varnothing$, which completes the proof
of (\ref{souscouple}).
\end{pf*}
\begin{Remark*}
On our way, we proved that for each $\lambda\in
\Lambda$,
\[
\lim_{t \to+\infty} \tilde{\P}_{\lambda}( \xi^{0}_{t} \ne
\varnothing,
\hat{\xi}^x_{t} \ne\varnothing, \xi^{0}_{t}\cap\hat{\xi
}^x_{t}=\varnothing)=0,\vadjust{\goodbreak}
\]
which is the essential ingredient in the proof of the complete
convergence Theorem~\ref{thCC}. One can refer to the article by
Durrett \cite{MR1117232} for the details in the case of the classical
contact process.
\end{Remark*}

%
\begin{appendix}\label{app}
\section*{\texorpdfstring{Appendix: Proof of almost subadditive
ergodic Theorems~\lowercase{\protect\ref{therg}} and \lowercase{\protect\ref{thergdeux}}}
{Appendix: Proof of almost subadditive ergodic Theorems 23 and 24}}

\begin{pf*}{Proof of Theorem \ref{therg}}
Let $a_p=C_p^{1/\alpha}$ and $u_n=\E[f_n]$; for every \mbox{$n,p \in\N$},
we have
$\E[r_{n,p}^+]\le(\E[(r_{n,p}^+)^{\alpha}])^{1/\alpha}\le
C_p^{1/\alpha}=a_p$,
hence,
\[
u_{n+p}\le u_n+u_p+\E[g_p]+\E[r_{n,p}]\le u_n+u_p+\E[g_p] +a_{p}.
\]
The general term of a convergent series tends to $0$, so
$C_p=o(p^{\alpha})$ or $a_p=o(p)$. Since $\frac{a_n+\E g_n}n$ tends to
$0$, the convergence of $u_n/n$ is classical (see Derriennic~\cite{MR704553}, e.g.).
The limit $\mu$ is finite because $u_n\ge
cn$ holds for each $n$.

We are going to show that $\underline f=\liminf_{n \to+\infty
}\frac{f_{n}}n $ stochastically dominates a~random variable whose mean
value is not less than $\mu$.

For every random variable $X$, let us denote by $\mathcal L(X)$ its law
under $\P$.
We denote by $\mathcal{K}$ the set of probability measures on $\R
_+^{\N
^*}$ whose marginals $m$ satisfy
\[
\forall t>0\qquad m(]t,+\infty[)\le\P(f_1+g_1>t/2)+C_1(2/t)^{\alpha}.
\]
Define, for $k \ge1$,
\[
\Delta_k=f_{k+1}-f_{k}
\]
and denote by $\Delta$ the process $\Delta=(\Delta_k)_{k \ge1}$.
For $k\in\N$, subadditivity ensures that $\Delta_k\le(f_1+g_1)\circ
\theta_k+r_{k,1}$,
hence, for each $t>0$,
\begin{eqnarray*}
\P(\Delta_k>t) & \le& \P\bigl((f_1+g_1)\circ\theta_k>t/2\bigr)+\P(r_{k,1}^+>t/2)
\\
& \le& \P(f_1+g_1>t/2)+C_1(2/t)^{\alpha}.
\end{eqnarray*}
This ensures that $\Delta\in\mathcal K$.

We denote by $s$ the shift operator $s((u_k)_{k \ge0})=(u_k)_{k \ge1}$,
and consider the sequence of probability measures on $\R^{\N^*}$
\[
(L_n)_{n \ge1}=\Biggl( \frac1n \sum_{j=1}^{n} \mathcal L(s^j \circ
\Delta) \Biggr)_{n \ge1}.
\]
Since $\mathcal{K}$ is convex and invariant by $s$, the sequence
$(L_n)_{n\ge1}$ is $\mathcal{K}$-valued.
Let \mbox{$n,k\ge1$}.
\begin{eqnarray*}
\int\pi_k(x) \,dL_{n}(x)
& = & \frac1{n} \sum_{j=1}^{n} \E\bigl(\pi_k(s^j \circ\Delta)\bigr) \\
& = & \frac1{n} \sum_{j=1}^{n} \E(f_{k+j+1}-f_{k+j})\\
& = &\frac1{n} (\E[f_{n+k+1}]-\E[f_{k+1}]).
\end{eqnarray*}
Let $M_k=\sup_{n\ge1} \frac1{n} |\E[f_{n+k+1}]-\E[f_{k+1}]|$.
The convergence of $u_n/n$ implies that $M_k$ is finite. Similarly, the
subadditivity gives
\begin{eqnarray*}
\int\pi_k^+(x) \,dL_{n}(x)
& = & \frac1{n} \sum_{j=1}^{n} \E\bigl(\pi_k^+(s^j \circ\Delta)\bigr) \\
& = & \frac1{n} \sum_{j=1}^{n} \E[(f_{k+j+1}-f_{k+j})^+]\\
&\le&\E[f_1^+]+\E[g_1]+a_1.
\end{eqnarray*}
Thus, we have
\begin{eqnarray*}\int|\pi_k(x)| \,dL_{n}(x)&\le& \int2\pi_k^+(x)\,
dL_{n}(x)+ \biggl|\int\pi_k(x) \,dL_{n}(x)\biggr|\\
&\le& M_k+2\E[f_1^+]+2\E
[g_1^+]+2a_1.
\end{eqnarray*}
Let $\mathcal{K}'$ be the family of laws $m$ on $\R^{\N^*}$ such that
for each $k$,
$\int|\pi_k| \,dm\le2 M_k+\E[f_1^+]+\E[g_1^+]+a_1$.
$\mathcal{K}'$ is compact for the topology of the convergence in law and
the sequence $(L_n)_{n\ge1}$ is $\mathcal{K}'$-valued. So, let
$\gamma
$ be a limit point of $(L_n)_{n \ge1}$ and $(n_k)_{k\ge1}$ a sequence
of indexes such that $L_{n_k}\Longrightarrow\gamma$.
By construction, $\gamma$ is invariant under the shift $s$.

Now, the sequence of the laws of the first coordinate $\pi_1(x)$ under
$(L_{n_k})_{k \ge0}$ weakly converges to the law of the first
coordinate under $\gamma$. Also, by definition of $\mathcal{K}$, the
positive parts of these elements form a uniformly integrable
collection, so
$\int\pi_1^+ \,d\gamma=\lim\int\pi_1^+ \,dL_{n_k}$.
However, the Fatou lemma tells us that
$\int\pi_1^- \,d\gamma\le\liminf_{k\to+\infty} \int\pi_1^-\,
dL_{n_k}$, hence, finally
\[
\int\pi_1 \,d\gamma\ge\liminf_{k\to+\infty}\int\pi_1\,
dL_{n_k}=\mu.
\]

Let\vspace*{1pt} $Y=(Y_k)_{k \ge1}$ be a process whose law is $\gamma$.
Since $\gamma$ is invariant under the shift~$s$, the Birkhoff theorem
tells us that
the sequence $(\frac1n \sum_{k=1}^{n}Y_k)_{n \ge1}$ a.s. converges to
a random variable $Y_{\infty}$, which then satisfies \mbox{$\E(Y_{\infty
})=\int\pi_1 \,d\gamma\ge\mu$}.

It remains to see that the law of $Y_{\infty}$ is stochastically
dominated by the law of $ \underline{f}=\liminf_{n \to+\infty
} \frac1{n}f_{n}$.
We will show that for each $a\in\R$, $\P(Y_{\infty}>a)\le\P
(\underline{f}>a)$.
By left-continuity, it is sufficient to prove the inequality in a dense
subset of $\R$. Thus, we can assume that $a$ is not an atom for the law
of $\underline{f}$:
\[
\{Y_{\infty}>a\}=\biggl\{\liminf_{n\to+\infty}\frac
{Y_1+\cdots+Y_n}{n}>a\biggr\}=\bigcup_{k\ge1}\biggl\{
\inf_{n\ge k}\frac{Y_1+\cdots+Y_n}{n}>a\biggr\}.
\]
Hence,
\begin{eqnarray*}
&&\P(Y_{\infty}>a) \\
&&\qquad = \limsup_{k\to+\infty}\P_{Y}\biggl(\inf_{n\ge
k} \frac{\pi_1+\cdots+\pi_n}{n}>a\biggr)\\
&&\qquad = \limsup_{k\to+\infty}\inf_{n\ge k}\P_{Y}
\biggl(\inf_{k\le i\le n} \frac{\pi_1+\cdots+\pi_i}{i}>a\biggr)\\
&&\qquad \le \limsup_{k\to+\infty}\inf_{n \ge k}
\liminf_{K \to+\infty}\frac{1}{n_K}\sum_{j=1}^{n_K}\P
\biggl(\inf_{k\le i\le n} \frac{\pi_1+\cdots+\pi_i}{i} \circ
s^j\circ\Delta>a\biggr).
\end{eqnarray*}
Let $\epsilon>0$. We have, for fixed $k,n,j$,
\begin{eqnarray*} & &\P\biggl(\inf_{k\le i\le n} \frac{\pi
_1+\cdots+\pi_i}{i} \circ s^j\circ\Delta>a\biggr)\\
&&\qquad = \P\biggl(\inf_{k\le i\le n} \frac{f_{i+j+1}-f_{j+1}}{i}
>a\biggr)\\
&&\qquad \le \P\biggl(\inf_{k\le i\le n} \frac{(f_{i}+g_i)\circ
\theta_{j+1}+r_{j+1,i}}{i} >a\biggr)\\
&&\qquad \le \P\biggl(\inf_{k\le i\le n} \frac{(f_{i}+g_i)\circ
\theta_{j+1}}{i} >a-\epsilon\biggr)+\P\biggl(\sup_{i\ge k}
\frac{r_{j+1,i}}{i} >\epsilon\biggr).
\end{eqnarray*}
On one hand, we have
\begin{eqnarray*}
\P\biggl(\sup_{i\ge k} \frac{r_{j+1,i}}{i} >\epsilon\biggr) &
\le&
\P\biggl(\sum_{i\ge k} \biggl(\frac{r^+_{j+1,i}}{i}
\biggr)^{\alpha} >\epsilon^{\alpha}\biggr)
\le\epsilon^{-\alpha}\sum_{i\ge k}\frac1{i^{\alpha}} \E
[(r^+_{j+1,i})^{\alpha}] \\
&\le& \epsilon^{-\alpha}\sum_{i\ge k} \frac{C_i}{i^{\alpha}}.
\end{eqnarray*}
We can note that this term does not depend on $j$ nor on $n$.
On the other hand,
\[
\P\biggl(\inf_{k\le i\le n} \frac
{(f_{i}+g_i)\circ\theta_{j+1}}{i} >a-\epsilon\biggr) = \P
\biggl(\inf_{k\le i\le n} \frac{f_{i}+g_i}{i} >a-\epsilon\biggr),
\]
which does not depend on $j$. Then, for each $\epsilon>0$, we have for
every $n,k$, with $n\ge k$,
\begin{eqnarray*}
&& \liminf_{K \to+\infty}\frac{1}{n_K}\sum_{j=1}^{n_K}\P
\biggl(\inf_{k\le i\le n} \frac{\pi_1+\cdots+\pi_i}{i} \circ
s^j\circ\Delta>a\biggr) \\
&&\qquad \le \epsilon^{-\alpha}\sum_{i\ge k} \frac{C_i}{i^{\alpha
}}+\P\biggl(\inf_{k\le i\le n} \frac{f_{i}+g_i}{i} >a-\epsilon
\biggr);
\end{eqnarray*}
next
\begin{eqnarray*}
&& \inf_{n\ge k} \liminf_{K \to+\infty}\frac
{1}{n_K}\sum_{j=1}^{n_K}\P
\biggl(\inf_{k\le i\le n} \frac{\pi_1+\cdots+\pi_i}{i} \circ
s^j\circ\Delta>a\biggr) \\
&&\qquad \le \epsilon^{-\alpha}\sum_{i\ge k} \frac{C_i}{i^{\alpha
}}+\inf_{n\ge k}\P\biggl(\inf_{k\le i\le n} \frac
{f_{i}+g_i}{i} >a-\epsilon\biggr).
\end{eqnarray*}
Finally,
\begin{eqnarray*}
\P(Y_{\infty}>a)
& \le& \limsup_{k\to+\infty}\inf_{n \ge k} \P
\biggl(\inf_{k\le i\le n} \frac{f_{i}+g_i}{i} >a-\epsilon
\biggr)+\limsup_{k\to+\infty}\epsilon^{-\alpha}\sum
_{i\ge k} \frac{C_i}{i^{\alpha}} \\
& \le& \limsup_{k\to+\infty} \P\biggl(\inf_{i\ge
k} \frac{f_{i}+g_i}{i} >a-\epsilon\biggr) \\
& \le& \P\biggl(\liminf_{i\to+\infty} \frac{f_{i}+g_i}{i}
>a-\epsilon\biggr)
= \P\biggl(\liminf_{i\to+\infty} \frac{f_{i}}{i} >a-\epsilon
\biggr),
\end{eqnarray*}
considering that $g_i/i$ almost surely converges to $0$.
Letting $\epsilon$ tend to zero, we obtain
\[
\P(Y_{\infty}>a)\le\P\biggl(\liminf_{i\to+\infty} \frac
{f_{i}}{i}\ge a \biggr)=\P(\underline{f}>a).
\]

It remains to see that $\underline f$ is invariant under the $\theta_n$'s.
Fix $n\ge1$. We have
\[
\E\Biggl[ \sum_{p=1}^{+\infty} \biggl(\frac{r^+_{n,p}}p
\biggr)^{\alpha}
\Biggr]
= \sum_{p=1}^{+\infty} \E\biggl[\biggl(\frac{r^+_{n,p}}p
\biggr)^{\alpha
} \biggr] \le\sum_{p=1}^{+\infty}\frac{C_p}{p^{\alpha}}<+\infty.
\]
Particularly, $\frac{r^+_{n,p}}p$ almost surely converges to $0$ when
$p$ tends to infinity.
Since $f_{n+p}\le f_n+f_p\circ\theta_n+g_p\circ\theta_n+r_{n,p}^+$,
dividing by $n+p$ and letting $p$ tend to $+\infty$, it comes that
\[
\underline{f} \le\underline{f} \circ\theta_n \qquad\mbox{a.s.}
\]
Since $\P$ is invariant under $\theta_n$, we classically conclude that
$\underline{f}$ is invariant under $\theta_n$.
\end{pf*}
\begin{Remark*}
In the present article, we made no use of the
possibility to take a nonzero $g_p$. In the case where the $(g_p)$ are
not zero, but the $r_{n,p}$'s are, we obtain a result which sounds a
bit like Theorem 3 in Sch\"{u}rger \cite{MR1127716}.
Like Sch\"{u}rger \cite{MR833959}, we use the idea of a coupling with a
stationarized process. This idea is due to Durrett \cite{MR586774} and
has been popularized by Liggett \cite{MR806224}. However, here there is
a refinement, because we directly establish a stochastic comparison
with the random variable $Y$, whereas previous papers establish a
stochastic comparisons with the whole process $(Y_n)_{n\ge1}$, that
admits $Y$ as its infimum limit.\vadjust{\goodbreak}

In the majority of almost subadditive ergodic theorems, almost sure
convergence requires strong conditions on the lack of subadditivity
(stationarity, e.g.).
Here we obtain an almost sure behavior by only considering a condition
on the moments (of order greater than 1) of the lack of subadditivity.
Besides, we know that bounding the first moment of the lack of
subadditivity is not sufficient to get an almost sure behavior (see the
remark by Derriennic \cite{MR704553} and the counter-example by
Derriennic and Hachem \cite{MR939537}).
\end{Remark*}

\begin{pf*}{Proof of Theorem \ref{thergdeux}}
It remains to prove that $ \E(\limsup_{n
\to+\infty} \frac{f_{n}}n ) \le\mu$.

We fix $k \ge1$. By subadditivity, we have for each $n\ge0$ and every
$0 \le r \le k-1$,
\begin{eqnarray*}
f_{nk+r}& \le& f_{nk}+(f_r+g_r)\circ\theta_{nk}+r_{nk,r}^+\\
& \le&\Biggl(\sum_{i=0}^{n-1}f_k\circ(\theta_k)^i\Biggr)+
\Biggl(f_{nk}-\sum_{i=0}^{n-1}f_k\circ(\theta_k)^i
\Biggr)^+\\
&&{}+(f_r+g_r)\circ
\theta_{nk}+r_{nk,r}^+.
\end{eqnarray*}
Since $\P$ is invariant under $\theta_k$, the Birkhoff theorem gives
the $L^1$ and almost-sure convergence
\[
\lim_{n \to+\infty} \frac1n \sum_{j=0}^{n-1} \frac{f_k \circ
(\theta
_k)^j}{k} =\frac{\E(f_k|\mathcal{I}_k)}k,
\]
where $\mathcal{I}_k$ is the $\sigma$-algebra of the $\theta
_k$-invariant events.
Let us now control the residual terms.
Since the finite collection $(f_r+g_r)_{0 \le r \le k-1}$ is
equi-integrable and $\P$ is invariant under $ \theta_k$, the collection
$(\sup_{0 \le r \le k-1}(f_r+g_r)\circ\theta_k^n)_{n \ge1}$
is equi-integrable, which ensures the almost sure and $L^1$ convergence
\[
\lim_{n \to+\infty} \frac{1}{n}\sup_{0 \le r \le k-1}
(f_r+g_r)\circ
(\theta_k)^n=0.
\]
We have $ \sum_{n=1}^{+\infty} \E[ (\frac
{r^+_{nk,r}}n)^{\alpha} ] \le\sum_{n=1}^{+\infty}\frac
{C_r}{n^{\alpha}}<+\infty$,
which implies, as previously, that $r_{nk,r}^+/n$ almost surely
converges to $0$.
Finally,
\[
\forall r\in\{0,\ldots,k-1\}\qquad \limsup_{n \to+\infty}\frac
{f_{nk+r}}{nk+r} \le\frac{\E[f_k|\mathcal{I}_k]}k,
\]
hence, $ \E[\limsup_{n \to+\infty}\frac
{f_{n}}{n}] \le\frac{\E[f_k]}k$.
We complete the proof by letting $k$ tend to $+\infty$.
\end{pf*}
\begin{Remark*}
When there is no lack of subadditivity, the
assumptions of Theorem \ref{thergdeux} obviously hold; thus we obtain a
subadditive ergodic theorem which sounds very much like Liggett's
\cite{MR806224}. However, these theorems are not strictly comparable,
in the
following sense that no one implies the other one.\vadjust{\goodbreak}

Indeed, extending a remark made by Kingman in his Saint-Flour's\break
course~\cite{MR0438477}, page 178, we can note that the assumption of
Kingman's original article [the stationarity of the doubly indexed
process $(X_{s,t})_{s\ge0,t\ge0}$] can be weakened in two different ways:
\begin{itemize}
\item Either assuming that for each $k$, the process
$(X_{(r-1)k,rk})_{r\ge1}$ is stationary; this assumption will be used
by Liggett \cite{MR806224}.
\item Or assuming that the law of $X_{n,n+p}$ does not depend on $p$.
That assumption, suggested by Hammersley and Welsh, is the one that we
use here, also used by Sch\"{u}rger in \cite{MR1127716}.
\end{itemize}
Note, however, that the special assumption of stationarity is used in
Liggett's proof \cite{MR806224} only in the so-called easy part, that
is, the bound for the supremum limit.

Kingman thought that the first set of assumptions surpassed the second
one, in view of possible applications. More than 30 years later, the
progresses of subadditive ergodic theorems, particularly about bounding
the infimum limit, lead to moderate this affirmation.
\end{Remark*}
\end{appendix}


%

\printaddresses

\end{document}